\documentclass[final,5p,times,twocolumn,authoryear]{elsarticle}
\usepackage{placeins}
\usepackage{amsmath}

\usepackage[bf]{caption}
\usepackage{longtable}
\usepackage{natbib}
\usepackage{rotating}
\usepackage{pdflscape}
\usepackage{appendix}
\usepackage{lipsum}
\usepackage{graphicx}
\usepackage{caption}
\usepackage{subcaption}

 \bibpunct[, ]{(}{)}{,}{a}{}{,}%

\usepackage{nccmath}
\usepackage{multirow}
\usepackage{array, booktabs}
\usepackage[hidelinks]{hyperref}
\usepackage{multicol}

\usepackage{tabularx}

\usepackage{amssymb}

\journal{}

\begin{document}

\begin{frontmatter}



\title{A computational study for the inventory routing problem}


\author[1]{Yasemin Malli}
\author[2]{Marco Laumanns}
\author[3]{Roberto Rossi}
\author[4]{Steven Prestwich}
\author[5]{S. Armagan Tarim}
 
\address[1]{Department of Management, Hacettepe University, Ankara, Turkey }
\address[2]{IBM Research-Zurich}
\address[3]{University of Edinburgh Business School, Edinburgh, UK}
\address[4]{Department of Computer Science, University College Cork, Ireland}
\address[5]{Cork University Business School, Cork, Ireland}
\begin{abstract}

In this work we compare several new computational approaches to an
inventory routing problem, in which a single product is shipped from a
warehouse to retailers via an uncapacitated vehicle.  We survey exact
algorithms for the Traveling Salesman Problem (TSP) and its
relaxations in the literature for the routing component.  For the
inventory control component we survey classical mixed integer linear
programming and shortest path formulations for inventory models.  We
present a numerical study comparing combinations of the two
empirically in terms of cost and solution time.

\end{abstract}

\begin{keyword}
inventory routing, traveling salesman problem, inventory control, mixed integer programming



\end{keyword}

\end{frontmatter}

\section{Introduction}

The inventory routing problem (IRP) is a distribution problem
combining inventory control and vehicle routing.  In an IRP products
are shipped from a warehouse to several geographically dispersed
retailers by means of vehicles (the routing component).  The warehouse
and retailers must manage their inventories without causing stock-outs
(the inventory control component).  IRPs try to answer the following
questions:
\begin{itemize}
\item
When should products be delivered to the retailers?
\item
What product quantities should be delivered to each visited retailer?
\item
Which routes should be taken?
\item
Which time periods should be replenishment periods?
\end{itemize}
The objective is to minimize the sum of expected inventory and
transportation costs during the planning horizon.

The routing component alone makes the IRP a challenging problem.  It
reduces to the traveling salesman problem (TSP) when the planning
horizon is one, the inventory costs are zero, the vehicle capacity is
infinite, and all customers need to be served.  Even with only one
customer some variants remain computationally hard 

The IRP has several variations: the planning horizon can be finite or
infinite; the number of products can be one or more; inventory holding
cost can be taken into account or not, and can be charged at the
warehouse, at the retailers or both; trucks, warehouses and retailers
may be capacitated or uncapacitated; and demands for each retailer can
be deterministic or stochastic, stationary or non-stationary.  The IRP
variation treated in this paper is as follows: there is a set of
retailers with deterministic but time varying demands; a single
product; trucks, warehouses and retailers have unlimited capacity;
dispatching a truck has a fixed cost; there are transportation costs
associated with the route followed by the truck to reach a given
retailer.  The problem is to determine the optimal replenishment plan
for each retailer and vehicle routes for each replenishment period,
while minimising total expected cost.

Most IRP research assumes capacitated warehouse, retailers and trucks,
but even without capacity constraints the problem is NP-hard because
of the TSP component.  It is surprising that no study in the
literature has investigated the uncapacitated case, which simply
integrates a TSP and a Wagner-Whitin inventory model.  There are
several formulations in the literature for both of these problems, and
investigating their possible combinations could yield useful guidance
for future IRP research.  Our contribution to the literature is as
follows:
\begin{itemize}
\item
We model the uncapacitated IRP routing component as a TSP and the
inventory component as a dynamic lot-sizing problem.
\item
We explore several combinations of TSP formulations and inventory
control models from the relevant literatures.
\item
We present a numerical study comparing these combinations empirically,
in terms of cost and solution time.
\end{itemize} 
The rest of paper is organized as follows.  Section 2 introduces the
relevant literature.  Section 3 explains the relevant TSP and
inventory control models.  Section 4 provides the definition of our
IRP.  Section 5 presents the numerical study.  Finally, Section 6
draws conclusions.

\section{Background}

The paper of \cite{federgruen1984} was among the first to investigate
the integration of inventory management and vehicle routing problems.
They consider a warehouse with several retailers and a scarce resource
available for a single period.  They use a fleet of vehicles with
limited capacity to deliver product.  The demand for each retailer is
assumed to be a random variable.  The objective is to determine
vehicle routes and replenishment quantities while minimizing the
expected sum of transportation and inventory costs, the latter
including holding and shortage costs.  They propose a non-linear mixed
integer programming model for this problem, and use generalized
Benders' decomposition that decomposes the problem into (i) an
inventory allocation problem with holding and shortage costs, and (ii)
a routing problem that can be computed as a TSP for each vehicle with
transportation costs.

\cite{dror1987} investigate reducing the long-term planning of
inventory routing problem to a short-term planning.  They consider the
case of capacitated vehicles, and present a mixed-integer linear
programming (MILP) model under deterministic demands.

\cite{campbell1998} explore how to model the long-term effects of
short-term decisions, and propose two short-term planning approaches.
They use a heuristic implemented in a rolling horizon framework that
determines a distribution strategy, while minimizing average
distribution cost over the planning horizon without causing stock-outs
at any retailer.

\cite{anily1990} analyze the case of a single product and an infinite
horizon.  The product is shipped from warehouse to retailers that have
deterministic demands by using a fleet of capacitated trucks.  They
propose a heuristic method to determine long-term replenishment and
routing plans while minimising transportation and inventory costs.

\cite{chandra1994} consider a plant that produces several products
delivered by a fleet of capacitated vehicles.  They present a
computational study which compares two approaches under deterministic
demands: one solves the production scheduling and routing problems
seperately, while the other combines the two problems into a single
model.

\cite{bertazzi2002} present a policy called the deterministic
order-up-to policy for an IRP.  The demand of each retailer is
deterministic and can vary over time.  Each retailer has a minimum and
maximum inventory level and must be visited before its inventory
reaches its minimum level.  When a retailer is visited the inventory
level is increased up to the maximum level.  They propose a heuristic
to determine a shipping policy that minimizes transportation cost and
inventory costs both at the supplier and the retailers.
\cite{archetti2007} present a MILP model for this problem and solve
the model optimally.

\cite{gaur2004} tackle a real-world problem for Albert Heijn (a
supermarket chain in the Netherlands).  They assume that a single
product with deterministic time-varying demand is delivered via a
fleet of capacitated vehicles.  They present a solution method for a
periodic IRP: first they determine the delivery times and vehicle
routes; then trucks are assigned to the routes; finally, workload
balancing is used to adjust the departure times.

\cite{bertazzi2008} analyse how transportation and inventory costs
affect the optimal cost, and the impact of vehicle capacity and
inventory holding capacity on distribution strategies, under
deterministic demands.

As can be seen from this survey, no previous studies have explored the
usefulness of TSP/Wagner-Whitin formulation combinations.

\section{Formulations of TSP and Inventory Control Models}

This section surveys known MILP formulations for the two components of
the IRP.

\subsection{Formulations of TSP}

A large number of TSP formulations have been presented in the
literature, and we survey several MILP formulations.

The objective is to find the shortest route that starts from an
initial city, visits every city exactly once and returns back to the
initial city.  The integer linear program (ILP) formulation of
\cite{dantzig1954} was one of the first.  $R
=\lbrace1,2,\ldots,r\rbrace$ is the set of cities which are visited
and $c_{ij}$ the distance between each pair of cities.  $x_{ij} $ is a
binary decision variable which takes the value 1 if and only if the
salesman travels from city $i$ to city $j$.  The formulation of
\cite{dantzig1954} is as follows:

\setlength{\belowdisplayskip}{2pt} \setlength{\belowdisplayshortskip}{2pt}
\setlength{\abovedisplayskip}{2pt} \setlength{\abovedisplayshortskip}{2pt}
\begin{flalign}\label{eq:eq1}
 \text{Minimize} \quad \quad  \sum_{i=1}^{r}\sum_{j=1}^{r}(x_{ij}c_{ij})\hfill  \\
\text{Subject to} \quad  \sum_{i = 1 i\neq j}^{r}x_{ij} = 1 \quad \forall  j \in R \label{eq:eq2}\hfill  \\
\quad \sum_{j = 1i\neq j}^{r}x_{ij} = 1 \quad \forall i \in R \label{eq:eq3}\hfill \\
\sum_{i,j\in S}x_{ij}\leq\vert S \vert-1 \quad \forall S\subseteq R, \quad 2 \leq S \leq r-2 \label{eq:eq4}\hfill \\
x_{ij} \in \lbrace0,1\rbrace\hfill
\end{flalign}

Eq.~\eqref{eq:eq2} and Eq.~\eqref{eq:eq3} are called assignment
constraints.  Eq.~\eqref{eq:eq2} ensures that the salesman must leave
each city once while Eq.~\eqref{eq:eq3} ensures that the salesman must
enter each city once.  Eq.~\eqref{eq:eq4} are called subtour
elimination constraints, and prevent the formation of subtours (tours
involving proper subsets of the cities).

\cite{miller1960} proposed a MILP formulation by introducing a new
continuous variable $u_{i}$ to reduce the number of subtour
elimination constraints.
$u_{i}$ denotes
sequence number in which retailer $i$ is visited.  While the objective
function Eq.~\eqref{eq:eq1} and assignment constraints
Eq.~\eqref{eq:eq2} and Eq.~\eqref{eq:eq3} remain the same, the subtour
eliminiation costraints are as follows:
 \begin{flalign}
u_{i}-u_{j}+rx_{ij}\leq r-1 \quad \forall i,j = 2,3,\ldots,r \quad i \neq j 
\end{flalign}

\cite{desrochers1991} proposed an alternative subtour elimination
constraints with a stronger relation between $u_{i}$ and $x_{ij}$
improving the constraints of \cite{miller1960}.
\begin{multline}
u_{i}-u_{j}+(r-1)x_{ij}-(r-3)x_{ji} \leq r-2 \\ \quad \forall i,j = 2,3,\ldots,r \quad i \neq j 
\end{multline}

\cite{gavish1978} presented a flow-based formulation which can be used
for some transportation problems as well as the TSP.  They introduced
new continuous variables $y_{ij}$ to represent the flow between cities
$i$ and $j$ to reformulate the subtour elimination constraints.  In
their formulation, called the single commodity flow formulation by
\cite{orman2006}, objective function Eq.~\eqref{eq:eq1} and assignment
constraints Eq.~\eqref{eq:eq2} and Eq.~\eqref{eq:eq3} are retained
while subtour elimination constraints as follows:
\begin{flalign}
y_{ij} \leq (r-1)x_{ij} \quad \forall  i \in R   \quad \forall  j \in R \quad i\neq j\\
\sum_{j \neq 1}^{r}y_{1j} = r-1\\
\sum_{j = 1i\neq j}^{r}y_{ij} - \sum_{j = 2i\neq j}^{r}y_{ji} = 1 \quad  i = 2,3,\ldots,r  \\
y_{ij} \geq 0
\end{flalign}

\cite{finke1984} proposed an integer formulation taking the form of
two commodity network flow problems.  
They introduced continuous variables $y_{ij}$ to represent the flow of
commodity 1 between cities $i$ and $j$, and$ z_{ij}$ to represent the
flow of commodity 2 between cities $i$ and $j$.  They retain the
objective Eq.~\eqref{eq:eq1} and assignment constraints
Eq.~\eqref{eq:eq2} and Eq.~\eqref{eq:eq3}, while the subtour
elimination constraints are:

\begin{flalign}
\sum_{j=2}^{r}(y_{1j}-y_{j1}) = (r-1)  \\
\sum_{j=2}^{r}(z_{1j}-z_{j1}) = -(r-1) \\
\sum_{j=1}^{r}(y_{ij}-y_{ji})= -1 \quad i = 2,\ldots,r \quad  \forall i \neq j  \\
\sum_{j=1}^{r}(z_{ij}-z_{ji}) = 1 \quad i = 2,\ldots,r \quad  \forall  \neq j  \\
\sum_{j=1}^{r}(y_{ij}+z_{ji}) = r-1 \quad i = 1,\ldots,r \quad  \forall i\neq j  \\
(y_{ij}+z_{ij}) = (r-1)x_{ij} \quad i,j = 1,\ldots,r \quad  \forall i\neq j  \\
 y_{ij}, z_{ij} \geq 0
\end{flalign}

\cite{bektacs2014} proposed a systematic way to find a generalization
of the inequalities of \cite{miller1960} (referred to as MTZ) and
\cite{desrochers1991} (referred to as DL), via a polyhedral approach
that studies and analyzes the convex hull of feasible sets.  Because
MTZ has a weak linear programming (LP) relaxation, they proposed new
inequalities with a strong LP relaxation.  First they implemented this
approach for two sets and reinvented the MTZ and DL inequalities.
Then they obtained a generalization of these inequalities which
implies the constraints proposed by \cite{dantzig1954}, referred to as
CLIQUE for three-node subsets.  A weaker version of constraints
\cite{dantzig1954} are the CIRCUIT inequalities which are given as
follows:
\begin{align}  \label{eq:eq5}
\sum x_{ij}\leq\vert C \vert -1 \quad \forall C\in G_{v_{1}}
\end{align}
where $G_{v_{1}}$ is the complete graph induced by the set of nodes
$v_{1}=2,\ldots,r$ and $C$ is any set of arcs defining an elementary
circuit in $G_{v_{1}}$.

They illustrate a two-node version of the CLIQUE constraints (2CLQ) as
follows:
\begin{align}
x_{ij}+x_{ji}\leq1 \quad i,j = 2,\ldots,r
\end{align}
The DL inequalities generalized to three nodes are:
\begin{multline} \label{eq:eq6}
u_{i}-u_{k}+(r-1)(x_{ij}+x_{jk})+(r-3)(x_{kj}+x_{ji})\\+rx_{ik}+(r-4)x_{ki}\leq2r-4
\end{multline}
\begin{multline}\label{eq:eq7}
2u_{i}-u_{j}-u_{k}+(2r-2)(x_{ij}+x_{ik})+(2r-8)(x_{ji}+x_{ki})+\\(2r-5)(x_{jk}+x_{kj})\leq 4r-10
\end{multline}
\begin{multline} \label{eq:eq8}
-2u_{i}+u_{j}+u_{k}+(2r-8)(x_{ij}+x_{ik})+(2r-2)(x_{ji}+x_{ki})+\\(2r-5)(x_{jk}+x_{kj})\leq 4r-10
\end{multline}
The CLIQUE constraints can be obtained for three nodes (3CLQ) by
adding Eq.~\eqref{eq:eq7} for triple $(i,j,k)$ to Eq.~\eqref{eq:eq8}
for the same triple.  The new constraints Eq.~\eqref{eq:eq6} and
Eq.~\eqref{eq:eq7} for triples $(i,j,k)$, $(k,i,j)$ and $(j,k,i)$ are
respectively called non-radical (NR) and radical (R) constraints.  A
different class of inequalities called 2PATH constraints are also
presented in \cite{bektacs2014}:
 
\begin{multline}  \label{eq:eq9}
u_{i}-u_{k}+(2r-3)x_{ik}+(r-4)x_{ki}+\\(r-1)(x_{ij}+x_{jk})\leq 2r-4 \quad i\neq j \neq k
\end{multline}
\begin{multline}  \label{eq:eq10}
u_{k}-u_{i}+(2r-7)x_{ik}+(r-1)x_{ki}+\\(r-4)(x_{ij}+x_{jk})\leq 2r-6  \quad i\neq j \neq k
\end{multline}

These constraints imply a lifted version of the CIRCUIT inequalities
Eq.~\eqref{eq:eq5}, and the summation of constraints
Eq.~\eqref{eq:eq9} and Eq.~\eqref{eq:eq10} for $(i,j,k)$ gives the
lifted circuit inequalities (L3).

\subsection{Inventory Control Models}

Inventory control has an important role in terms of matching supply
and demand effectively.  Its aim is to determine the time and quantity
of replenishment orders while minimizing inventory costs.  In this
section we focus on two deterministic inventory models for the
well-known \cite{wagner1958} problem.

We first give a MILP formulation for the \cite{wagner1958} problem.
Consider an $N$-period planning horizon indexed by $t \in
\lbrace1,2,\ldots,N\rbrace$.  Assume that demand $d_{t}$ in period $t$
is deterministic.  A fixed holding cost $h$ is incurred if any item is
carried from one period to the next.  A unit cost $v$ is incurred for
each item ordered.  A fixed ordering cost $K$ is incurred each time an
order is placed.  The decision variables in this problem are as
follows: $q_{t}$ represents the quantity ordered at period $t$.
$I_{t}$ represents inventory level at period $t$.  $\delta_{t}$ is a
binary decision variable that takes value 1 if and only if an order is
placed at period $t$.  The formulation is:

\begin{flalign}
\text{Minimize} \quad \quad  \sum_{t=1}^{n}(K\delta_{t}+hI_{t}+vq_{t})\end{flalign}
\begin{equation}
\text{Subject to}\quad \quad   I_{t}=\begin{cases}
    I_{0}+q_{t}-d_{t} , & \text{if $t = 1$}.\\
    I_{t-1}+q_{t}-d_{t}, & \text{otherwise}.  
  \end{cases}\\
  \quad t = 1,2,\ldots,n \\
 \end{equation}
\begin{flalign} 
q_{t} \leq \sum_{k=t}^{n}d_{k}\delta_{k}\quad \textrm{for} \quad t = 1,2,\ldots,n  \\
I_{t},q_{t} \geq 0, \quad t = 1,2,\ldots,n 
\end{flalign}

The dynamic programming formulation proposed by \cite{wagner1958} can
also be formulated as shortest path network flow.  The objective of
the shortest path problem is to determine the shortest route between
source and destination.  The network is a directed acyclic graph with several nodes between source and
destination.  It consists of $N+1$ nodes representing periods to take
into account the cost of the last period.  The first period
represented by node 1 is the source and the destination is node $N+1$.
Arc $(i,j)$ is traversed if an order is placed at period $i$ to cover
demands from period $i$ to period $j-1$, and at period $j$ the next
order is placed.  So $\omega(i,j)$ can refer to the inventory cost
including ordering costs and holding cost incurred from period $i$ to
period $j-1$ (see \cite{muckstadt2010}).

\begin{align}
\omega(i,j) = K + \sum_{k=i}^{j}vd_{k}+h\sum_{k=i}^{j}(k-i)d_{k} \quad i,j = 1,\ldots,n+1 \quad i < j
\end{align}

The decision variable $w_{ij}$ determines whether an order is placed
at period $i$ to satisfy demand of period $i,i+1,\ldots,j-1$.  Hence
the shortest path LP formulation for the inventory control problems is
given as follows:
 
Minimize
\begin{align}
\sum_{i,j\in N}\omega(i,j)w_{ij}
\end{align}
Subject to
\begin{equation}
\sum_{j\in N+1}w_{ij} - \sum_{j\in N+1}w_{ji} = \begin{cases}
 1  & \text{if $i = 1$}\\
 -1  & \text{if $i = n+1$}\\
  0  & \text{Otherwise}\\  \end{cases}
\quad \quad i,j = 1, \ldots, n+1 \\
\end{equation}

\begin{align}
w_{ij} \geq 0
\end{align}

\section{Definition of the Inventory Routing Problem}

We consider a system which a single product is shipped from a
warehouse (W) to a set of retailers $R = \lbrace1,2,\ldots,r\rbrace$
over a finite planning horizon $N$.  We assume that for each period $t
\in \lbrace1,2,\ldots,N\rbrace$ the demand $d_{ti}$ of any retailer $i
\in \lbrace1,2,\ldots,R\rbrace$ is known.  A truck
departs from the warehouse, visits all the retailers which are
replenished once, then returns back to the warehouse via the route
with minimum distance.  The truck does not necessarily visit all
retailers at each period.  It is assumed that the capacities of the
warehouse, truck, and retailers are unlimited.

The parameters of the problem are defined as follows.  A holding cost
$h_{i}$ is incurred if any item is carried from one period to the next
by any retailer (we ignore unit ordering cost $v$ in this problem).  A
fixed ordering cost $K$ is incurred each time an order is placed by
any retailer.  The distance $c_{ij}$ between retailers $i$ and $j$ is
known.  The decision variables of the problem are as follows:
$x_{tij}$ is a binary decision variable from the TSP component that
determines whether the truck travels from city $i$ to $j$ in period
$t$.  From the inventory component, $q_{ti}$ represents the quantity
ordered at retailer $i$ in period $t$, and $I_{ti}$ represents the
inventory level of retailer $i$ at the end of period $t$.  Retailers
have zero initial and final inventory levels. Also from the inventory component,
$\delta_{ti}$ is a binary decision variable which takes the value 1 if
an order is placed by retailer $i$ in period $t$, and $w_{itk}$
determines whether an order is placed at period $t$ to satisfy demand
at period $t,t+1,\ldots,k-1$ for retailer $i$ in the shortest path
formulation.

We aim to find the optimal replenishment plan and vehicle routes that
minimize the total cost, including inventory cost and transportation
cost, without stock-outs at any retailers.  We use TSP formulations to
determine routing for each replenishment period and inventory control
models to determine replenishment quantity for each retailer mentioned
in Section 3.  We combine the MTZ, DL and TSP models proposed by
\cite{bektacs2014} with the \citeauthor{wagner1958}'s MILP inventory model by using
assignment constraints (Eq.~\eqref{eq:eq2} and Eq.~\eqref{eq:eq3}).
\begin{flalign}
\sum_{i = 1 i\neq j}^{r+1}x_{ij} = \delta_{tj} \quad \forall  j \in R\cup \lbrace W \rbrace \label{eq:CMILP+eq2}\hfill  \\
\quad \sum_{j = 1i\neq j}^{r+1}x_{ij} = \delta_{ti} \quad \forall i \in R\cup \lbrace W \rbrace \label{eq:CMILP+eq3}\hfill 
\end{flalign}
Eq.~\eqref{eq:CMILP+eq2} and Eq.~\eqref{eq:CMILP+eq3} ensure that if
there is an order at a retailer then the truck must visit that
retailer.  For the combination of single commodity (SC) flow
formulation and \citeauthor{wagner1958}'s MILP,  we use the
following equations in addition to assignment constraints:
\begin{flalign}
\sum_{j = 1i\neq j}^{r+1}y_{ij} - \sum_{j = 2i\neq j}^{r+1}y_{ji} = \delta_{ti} \quad  i = 2,3,\ldots,r+1 \\
\sum_{j \neq 1}^{r}y_{1j} = (r-1)\delta_{tj}
\end{flalign}
For the combination of the two commodity (2C) flow formulation and
\citeauthor{wagner1958}'s MILP, the assignment constraints remain the same while the
linking constraints are:
\begin{flalign}
\sum_{j=2}^{r+1}(y_{1j}-y_{j1}) = (r-1)\delta_{ti} \\
\sum_{j=2}^{r+1}(z_{1j}-z_{j1}) = -(r-1) \delta_{ti}\\
\sum_{j=1}^{r+1}(y_{ij}-y_{ji})= -1\delta_{ti} \quad i = 2,\ldots,r+1 \quad  \forall i \neq j  \\
\sum_{j=1}^{r+1}(z_{ij}-z_{ji}) = 1\delta_{ti} \quad i = 2,\ldots,r+1 \quad  \forall  \neq j  \\
\sum_{j=1}^{r+1}(y_{ij}+z_{ji}) = (r-1)\delta_{ti} \quad i = 1,\ldots,r+1 \quad  \forall i\neq j 
\end{flalign}

For models relying on the shortest path network flow formulation of \cite{wagner1958} model the linking constraints are
redefined by using decision variables $w_{itk}$.

\section{Numerical Study}

In this section we conduct two types of numerical study to compare the
solution times of optimal algorithms, and evaluate the strengths of
linear relaxation models.

First we consider the instance set proposed by \cite{archetti2007}.
There are two types of planning horizon lengths $N \in \lbrace 3,6
\rbrace$.  We have one warehouse and ten possible values of retailers
such that $r=5k$ ($k = 1,2,\ldots,10$) when the planning horizon is
$3$.  If the planning horizon is $6$ then we have six possible values
of retailers such that $r=5k$ ($k = 1,2,\ldots,6$).  We assume that
the warehouse has sufficient product for its retailers, and that there
is no initial inventory at any retailer.  The demand for each retailer
is deterministic and constant for each period.  The demand is a
uniformly distributed random integer in the interval $[10,100]$.
There are different holding costs for each retailer, uniformly
randomly distributed in the interval $[0.01,0.05]$.  There is no
holding cost for the warehouse.  We consider inventory costs at
retailers and transportation costs between warehouse and retailers.
The fixed ordering cost is $10$.  The distance $c_{ij}$ between
retailers $i$ and $j$ is given by Euclidean geometry:
$c_{ij}=\sqrt{(x_{i}-x_{j})^{2}+(y_{i}-y_{j})^{2}}$ where coordinates
$(x_{i},y_{i})$ and $(x_{j},y_{j})$ are uniformly distributed random
integers in the interval $[0,500]$.  In this numerical design the
complexity of TSP dominates that of the inventory component, and the
computational times depend mainly on the configuration of retailers
(Tables \ref{table:tab1}---\ref{table:bektasSP1}).

In our second numerical design, based on instances derived from the
academic literature on the TSP and inventory management, we focus
instead on problems in which inventory control occurs over a longer
planning horizon of 15.  We have one warehouse and one type of
retailer set of 16 retailers \citep{tsplib}.  A truck with unlimited
capacity delivers products.  We assume that the demand of any retailer
is deterministic for each period, but with three demand patterns.  In
the first demand pattern each retailer has the same stationary demand
(STA) of 100 in each period.  In the second demand pattern retailers
have different stationary demand patterns: the demand for five
retailers is 100 for each period, the demand for six other retailers
is 50 and for the rest 75 for each period.  The last demand pattern is
a combination of stationary demand pattern, two life cycle patterns
(LCY1 and LCY2), two sinusoidal patterns (SIN1 and SIN2) and a random
pattern (RAND).  We adopt demand patterns from
\cite{rossi2015piecewise}. Figure \ref{fig:graph} illustrates these demand
patterns.  We consider three types of fixed ordering cost $K$: first,
no fixed ordering cost; secondly, a fixed ordering cost of 1000 for
each retailer; finally, a fixed ordering cost of 500 for five
retailers, 1000 for six other retailers and 2000 for the remaining
five retailers.  In this numerical design we also have dispatching
cost $D$ which is incurred if the truck visits any nonzero number of
retailers.  We have two types of dispatching cost $D \in \lbrace
0,15000 \rbrace$.  We assume that the warehouse has sufficient product
to deliver to retailers, and that there is no initial inventory at any
retailer.  There is no holding cost for the warehouse, and the holding
cost is 1 for each retailer.  For each combination of $K$, $D$ and
demand patterns we generate 18 scenarios.

The results are presented in Tables
\ref{table:CMILP}--\ref{table:SPbektas} which show the computational
times (in seconds) and optimality gap (\%) of different combinations.
We set a time limit of 1 hour and record the optimality gap if a
problem is not solved optimally.  CMILP and SP respectively denote
classical MILP and shortest path formulations for inventory control.
MTZ, MTZ+2Clq, DL, SC and 2C denote the exact TSP algorithms mentioned
in Section 3.  In Scenarios 1--3 we assume that there is no fixed
ordering cost and dispatching cost.  Retailers have three types of
demand pattern; same STA demand, different STA demand and combination
of STA, LCY1, LCY2, SIN1, SIN2 and RAND.  In Scenarios 4--6 there is a
fixed ordering cost of 1000 and no dispatching cost.  Three demand
patterns are used.  In Scenarios 7--9 we have the same demand patterns
and retailers have different fixed ordering costs: 500 for five
retailers, 1000 for six others and 2000 for the remaining five
retailers.  There is no dispatching cost.  In Scenarios 10--18 we
assume that the dispatching cost is 15000 and we use the same
combinations as in Scenarios 1--9.

Table \ref{table:CMILP} and Table \ref{table:SP} show the same
information.  Table \ref{table:CMILP} shows that the CMILP+MTZ,
CMILP+MTZ+2CLQ and CMILP+DL combinations do not yield optimal
solutions within the time limit.  However, the results of Scenarios
4--9 (in the presence of fixed ordering cost only) are very close to
optimal.  CMILP+SC yields optimal solutions for all scenarios except
1.  CMILP+2C is very close to optimal in most scenarios.  Table
\ref{table:SP} provides similar results but its computational times
are longer.

Tables \ref{table:LPCMILP} and \ref{table:LPSP} show the optimality
gaps for the LP relaxations of the combination of inventory models and
TSP algorithms.  As expected, the LP relaxation of CMILP+MTZ is
furthest from optimal.  The LP relaxations of CMILP+MTZ+2CLQ and
CMILP+DL combinations have the same optimality gap.  As indicated in
\cite{bektacs2014} the LP relaxation of CMILP+MTZ+2CLQ provides a
better approximation than that of CMILP+MTZ.  The LP relaxations of
CMILP+SC and CMILP+2C combinations yield the same optimality gap.
Although the LP relaxations of CMILP+SC and CMILP+2C are slightly
worse than those of CMILP+MTZ+2CLQ and CMILP+DL for Scenarios 1--9,
the LP relaxations of CMILP+SC and CMILP+2C are best for Scenarios
10--18.  Regarding computation times, the LP relaxations of
CMILP+MTZ+2CLQ and CMILP+DL are faster than those of CMILP+SC and
CMILP+2C for Scenarios 1--12.  We also observe that the computation
times for the LP relaxation of CMILP+MTZ are greater with the third
demand pattern and fixed ordering cost (Scenarios 6, 9, 15 and 18).
If we analyse the LP relaxations, SP+SC and SP+2C
LP relaxations yield the best approximations for all scenarios, while 
SP+MTZ gives the biggest optimality gap.  The computation times are similar for
all combinations and scenarios.

Tables \ref{table:gapcmılp} and \ref{table:gapSP} show the optimality
gaps for the LP relaxations of the inventory models and TSP
formulations proposed by \cite{bektacs2014}.  Tables
\ref{table:CMILPbektas} and \ref{table:SPbektas} show computation
times for these combinations.  From Table \ref{table:gapcmılp} and
Table \ref{table:gapSP} we see that most combinations yield the same
optimality gap.  While the LP relaxation of CMILP+TSP formulations
proposed by \cite{bektacs2014} provide better approximations than the
LP relaxation of CMILP+MTZ, CMILP+MTZ+2CLQ and CMILP+DL, the
computational times for the LP relaxations of the CMILP+TSP
formulations proposed by \cite{bektacs2014} are very high.  The LP
relaxation of the SP+TSP formulations proposed by \cite{bektacs2014}
give slightly better approximations than those of SP+MTZ, SP+MTZ+2CLQ
and SP+DL.  Comparing the combinations of CMILP and SP with TSP
formulations proposed by \cite{bektacs2014}, the CMILP combinations
give better approximations but greater computational times.

\section{Conclusion}

In this paper, we presented combinations of optimal and near-optimal
TSP algorithms and inventory models from the literature, to solve an
inventory routing problem under deterministic demand over a finite
planning horizon.  

We tested these algorithms empirically on two types
of scenario and found the following results:
\begin{itemize}
\item
For the combination of exact TSP algorithms and inventory models,
CMILP+TSP combinations are faster than SP+TSP combinations and the
fastest algorithm is CMILP+SC.
\item
CMILP+SC and SP+SC combinations are more sensitive as fixed ordering
and dispatching cost increase.
\item
The combinations of TSP formulations proposed by \cite{bektacs2014}
with inventory models provide better approximations than MTZ
formulations.
\end{itemize}   

\bibliography{references}

\begin{thebibliography}{20}
\expandafter\ifx\csname natexlab\endcsname\relax\def\natexlab#1{#1}\fi
\providecommand{\url}[1]{\texttt{#1}}
\providecommand{\href}[2]{#2}
\providecommand{\path}[1]{#1}
\providecommand{\DOIprefix}{doi:}
\providecommand{\ArXivprefix}{arXiv:}
\providecommand{\URLprefix}{URL: }
\providecommand{\Pubmedprefix}{pmid:}
\providecommand{\doi}[1]{\href{http://dx.doi.org/#1}{\path{#1}}}
\providecommand{\Pubmed}[1]{\href{pmid:#1}{\path{#1}}}
\providecommand{\bibinfo}[2]{#2}
\ifx\xfnm\relax \def\xfnm[#1]{\unskip,\space#1}\fi
\bibitem[{Anily and Federgruen(1990)}]{anily1990}
\bibinfo{author}{Anily, S.}, \bibinfo{author}{Federgruen, A.},
  \bibinfo{year}{1990}.
\newblock \bibinfo{title}{One warehouse multiple retailer systems with vehicle
  routing costs}.
\newblock \bibinfo{journal}{Management Science} \bibinfo{volume}{36},
  \bibinfo{pages}{92--114}.
\bibitem[{Archetti et~al.(2007)Archetti, Bertazzi, Laporte and
  Speranza}]{archetti2007}
\bibinfo{author}{Archetti, C.}, \bibinfo{author}{Bertazzi, L.},
  \bibinfo{author}{Laporte, G.}, \bibinfo{author}{Speranza, M.G.},
  \bibinfo{year}{2007}.
\newblock \bibinfo{title}{A branch-and-cut algorithm for a vendor-managed
  inventory-routing problem}.
\newblock \bibinfo{journal}{Transportation Science} \bibinfo{volume}{41},
  \bibinfo{pages}{382--391}.
\bibitem[{Bekta{\c{s}} and Gouveia(2014)}]{bektacs2014}
\bibinfo{author}{Bekta{\c{s}}, T.}, \bibinfo{author}{Gouveia, L.},
  \bibinfo{year}{2014}.
\newblock \bibinfo{title}{Requiem for the miller--tucker--zemlin subtour
  elimination constraints?}
\newblock \bibinfo{journal}{European Journal of Operational Research}
  \bibinfo{volume}{236}, \bibinfo{pages}{820--832}.
\bibitem[{Bertazzi et~al.(2002)Bertazzi, Paletta and Speranza}]{bertazzi2002}
\bibinfo{author}{Bertazzi, L.}, \bibinfo{author}{Paletta, G.},
  \bibinfo{author}{Speranza, M.G.}, \bibinfo{year}{2002}.
\newblock \bibinfo{title}{Deterministic order-up-to level policies in an
  inventory routing problem}.
\newblock \bibinfo{journal}{Transportation Science} \bibinfo{volume}{36},
  \bibinfo{pages}{119--132}.
\bibitem[{Bertazzi et~al.(2008)Bertazzi, Savelsbergh and
  Speranza}]{bertazzi2008}
\bibinfo{author}{Bertazzi, L.}, \bibinfo{author}{Savelsbergh, M.},
  \bibinfo{author}{Speranza, M.G.}, \bibinfo{year}{2008}.
\newblock \bibinfo{title}{Inventory routing}, in: \bibinfo{booktitle}{The
  vehicle routing problem: latest advances and new challenges}.
  \bibinfo{publisher}{Springer}, pp. \bibinfo{pages}{49--72}.
\bibitem[{Campbell et~al.(1998)Campbell, Clarke, Kleywegt and
  Savelsbergh}]{campbell1998}
\bibinfo{author}{Campbell, A.}, \bibinfo{author}{Clarke, L.},
  \bibinfo{author}{Kleywegt, A.}, \bibinfo{author}{Savelsbergh, M.},
  \bibinfo{year}{1998}.
\newblock \bibinfo{title}{The inventory routing problem}, in:
  \bibinfo{booktitle}{Fleet management and logistics}.
  \bibinfo{publisher}{Springer}, pp. \bibinfo{pages}{95--113}.
\bibitem[{Chandra and Fisher(1994)}]{chandra1994}
\bibinfo{author}{Chandra, P.}, \bibinfo{author}{Fisher, M.L.},
  \bibinfo{year}{1994}.
\newblock \bibinfo{title}{Coordination of production and distribution
  planning}.
\newblock \bibinfo{journal}{European Journal of Operational Research}
  \bibinfo{volume}{72}, \bibinfo{pages}{503--517}.
\bibitem[{Dantzig et~al.(1954)Dantzig, Fulkerson and Johnson}]{dantzig1954}
\bibinfo{author}{Dantzig, G.}, \bibinfo{author}{Fulkerson, R.},
  \bibinfo{author}{Johnson, S.}, \bibinfo{year}{1954}.
\newblock \bibinfo{title}{Solution of a large-scale traveling-salesman
  problem}.
\newblock \bibinfo{journal}{Journal of the operations research society of
  America} \bibinfo{volume}{2}, \bibinfo{pages}{393--410}.
\bibitem[{Desrochers and Laporte(1991)}]{desrochers1991}
\bibinfo{author}{Desrochers, M.}, \bibinfo{author}{Laporte, G.},
  \bibinfo{year}{1991}.
\newblock \bibinfo{title}{Improvements and extensions to the
  miller-tucker-zemlin subtour elimination constraints}.
\newblock \bibinfo{journal}{Operations Research Letters} \bibinfo{volume}{10},
  \bibinfo{pages}{27--36}.
\bibitem[{Dror and Ball(1987)}]{dror1987}
\bibinfo{author}{Dror, M.}, \bibinfo{author}{Ball, M.}, \bibinfo{year}{1987}.
\newblock \bibinfo{title}{Inventory /routing: Reduction from an annual to a
  short-period problem}.
\newblock \bibinfo{journal}{Naval Research Logistics (NRL)}
  \bibinfo{volume}{34}, \bibinfo{pages}{891--905}.
\bibitem[{Federgruen and Zipkin(1984)}]{federgruen1984}
\bibinfo{author}{Federgruen, A.}, \bibinfo{author}{Zipkin, P.},
  \bibinfo{year}{1984}.
\newblock \bibinfo{title}{A combined vehicle routing and inventory allocation
  problem}.
\newblock \bibinfo{journal}{Operations Research} \bibinfo{volume}{32},
  \bibinfo{pages}{1019--1037}.
\bibitem[{Finke et~al.(1984)Finke, Claus and Gunn}]{finke1984}
\bibinfo{author}{Finke, G.}, \bibinfo{author}{Claus, A.},
  \bibinfo{author}{Gunn, E.}, \bibinfo{year}{1984}.
\newblock \bibinfo{title}{A two-commodity network flow approach to the
  traveling salesman problem}.
\newblock \bibinfo{journal}{Congressus Numerantium} \bibinfo{volume}{41},
  \bibinfo{pages}{167--178}.
\bibitem[{Gaur and Fisher(2004)}]{gaur2004}
\bibinfo{author}{Gaur, V.}, \bibinfo{author}{Fisher, M.L.},
  \bibinfo{year}{2004}.
\newblock \bibinfo{title}{A periodic inventory routing problem at a supermarket
  chain}.
\newblock \bibinfo{journal}{Operations Research} \bibinfo{volume}{52},
  \bibinfo{pages}{813--822}.
\bibitem[{Gavish and Graves(1978)}]{gavish1978}
\bibinfo{author}{Gavish, B.}, \bibinfo{author}{Graves, S.C.},
  \bibinfo{year}{1978}.
\newblock \bibinfo{title}{The travelling salesman problem and related
  problems}.
\newblock \bibinfo{type}{Technical Report}. Massachusetts Institute of
  Technology, Operations Research Center.
\bibitem[{Miller et~al.(1960)Miller, Tucker and Zemlin}]{miller1960}
\bibinfo{author}{Miller, C.E.}, \bibinfo{author}{Tucker, A.W.},
  \bibinfo{author}{Zemlin, R.A.}, \bibinfo{year}{1960}.
\newblock \bibinfo{title}{Integer programming formulation of traveling salesman
  problems}.
\newblock \bibinfo{journal}{Journal of the ACM (JACM)} \bibinfo{volume}{7},
  \bibinfo{pages}{326--329}.
\bibitem[{Muckstadt and Sapra(2010)}]{muckstadt2010}
\bibinfo{author}{Muckstadt, J.A.}, \bibinfo{author}{Sapra, A.},
  \bibinfo{year}{2010}.
\newblock \bibinfo{title}{Principles of inventory management: When you are down
  to four, order more}.
\newblock \bibinfo{publisher}{Springer Science \& Business Media}.
\bibitem[{Orman and Williams(2006)}]{orman2006}
\bibinfo{author}{Orman, A.}, \bibinfo{author}{Williams, H.P.},
  \bibinfo{year}{2006}.
\newblock \bibinfo{title}{A survey of different integer programming
  formulations of the travelling salesman problem}.
\newblock \bibinfo{journal}{Optimisation, economics and financial analysis.
  Advances in computational management science} \bibinfo{volume}{9},
  \bibinfo{pages}{93--106}.
\bibitem[{Rossi et~al.(2015)Rossi, Kilic and Tarim}]{rossi2015piecewise}
\bibinfo{author}{Rossi, R.}, \bibinfo{author}{Kilic, O.A.},
  \bibinfo{author}{Tarim, S.A.}, \bibinfo{year}{2015}.
\newblock \bibinfo{title}{Piecewise linear approximations for the
  static--dynamic uncertainty strategy in stochastic lot-sizing}.
\newblock \bibinfo{journal}{Omega} \bibinfo{volume}{50},
  \bibinfo{pages}{126--140}.
\bibitem[{TSPLIB()}]{tsplib}
TSPLIB, \bibinfo{year}{1997}.
\newblock \bibinfo{title}{{TSPLIB}}.
\newblock
  \bibinfo{note}{Http://elib.zib.de/pub/mp-testdata/tsp/tsplib/tsplib.html}.
\bibitem[{Wagner and Whitin(1958)}]{wagner1958}
\bibinfo{author}{Wagner, H.M.}, \bibinfo{author}{Whitin, T.M.},
  \bibinfo{year}{1958}.
\newblock \bibinfo{title}{Dynamic version of the economic lot size model}.
\newblock \bibinfo{journal}{Management Science} \bibinfo{volume}{5},
  \bibinfo{pages}{89--96}.

\end{thebibliography}
\bibliographystyle{elsarticle-harv}

\renewcommand\appendix{\par
\setcounter{section}{0}%
\setcounter{subsection}{0}%
\setcounter{table}{0}
\setcounter{figure}{0}
\gdef\thetable{\Alph{table}}
\gdef\thefigure{\Alph{figure}}

\gdef\thesection{\Alph{section}}
\setcounter{section}{1}}
\appendix

\section*{Appendix}

In this appendix we illustrate demand patterns used in our computational study and the results of our experiments.

\begin{footnotesize}

\begin{figure*}[h!]
        \centering
        
        \begin{subfigure}[b]{0.45\textwidth}
        
                \includegraphics[width=\textwidth]{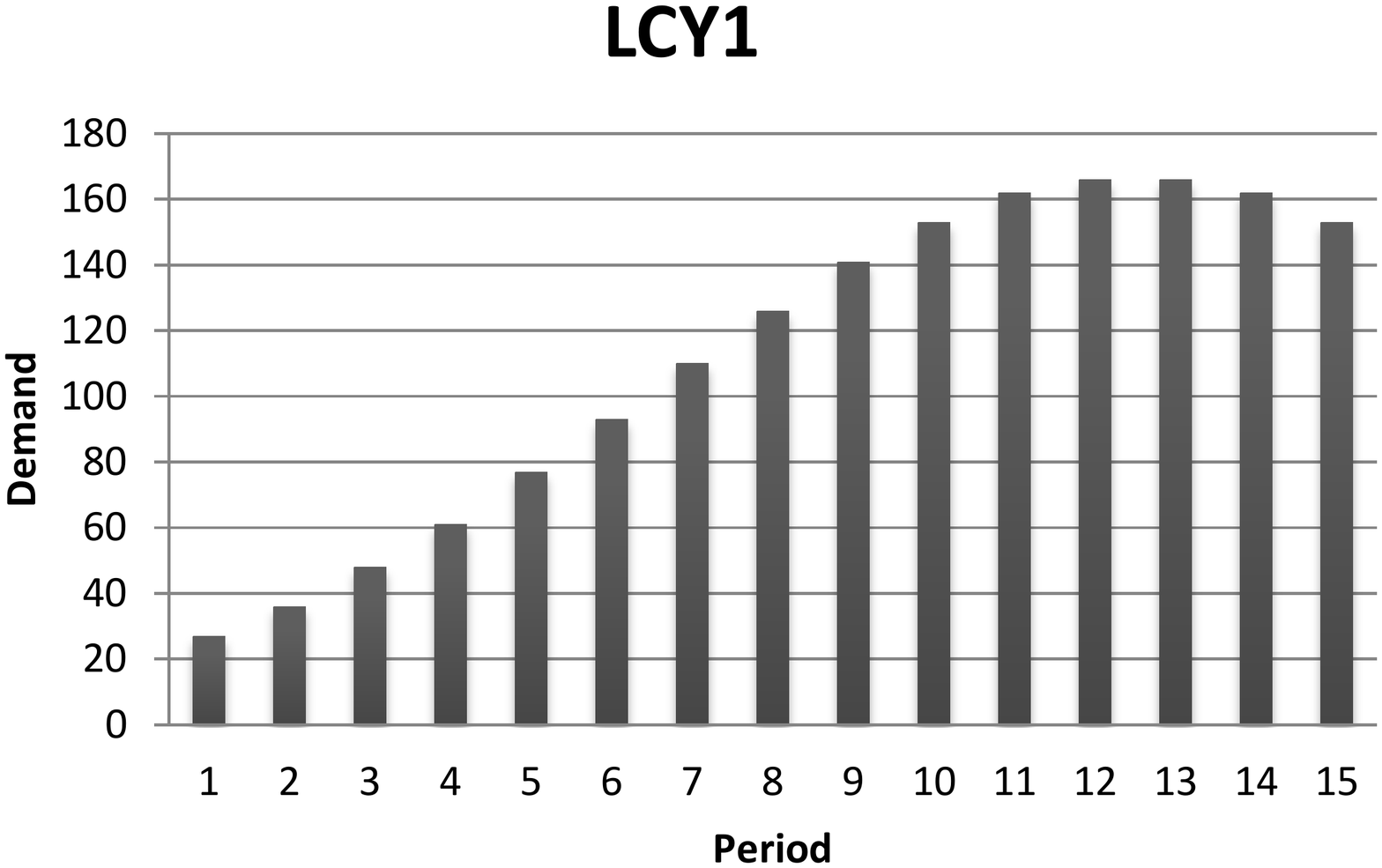}
                
        \end{subfigure}%
        ~ 
        \begin{subfigure}[b]{0.45\textwidth}
       
                \includegraphics[width=\textwidth]{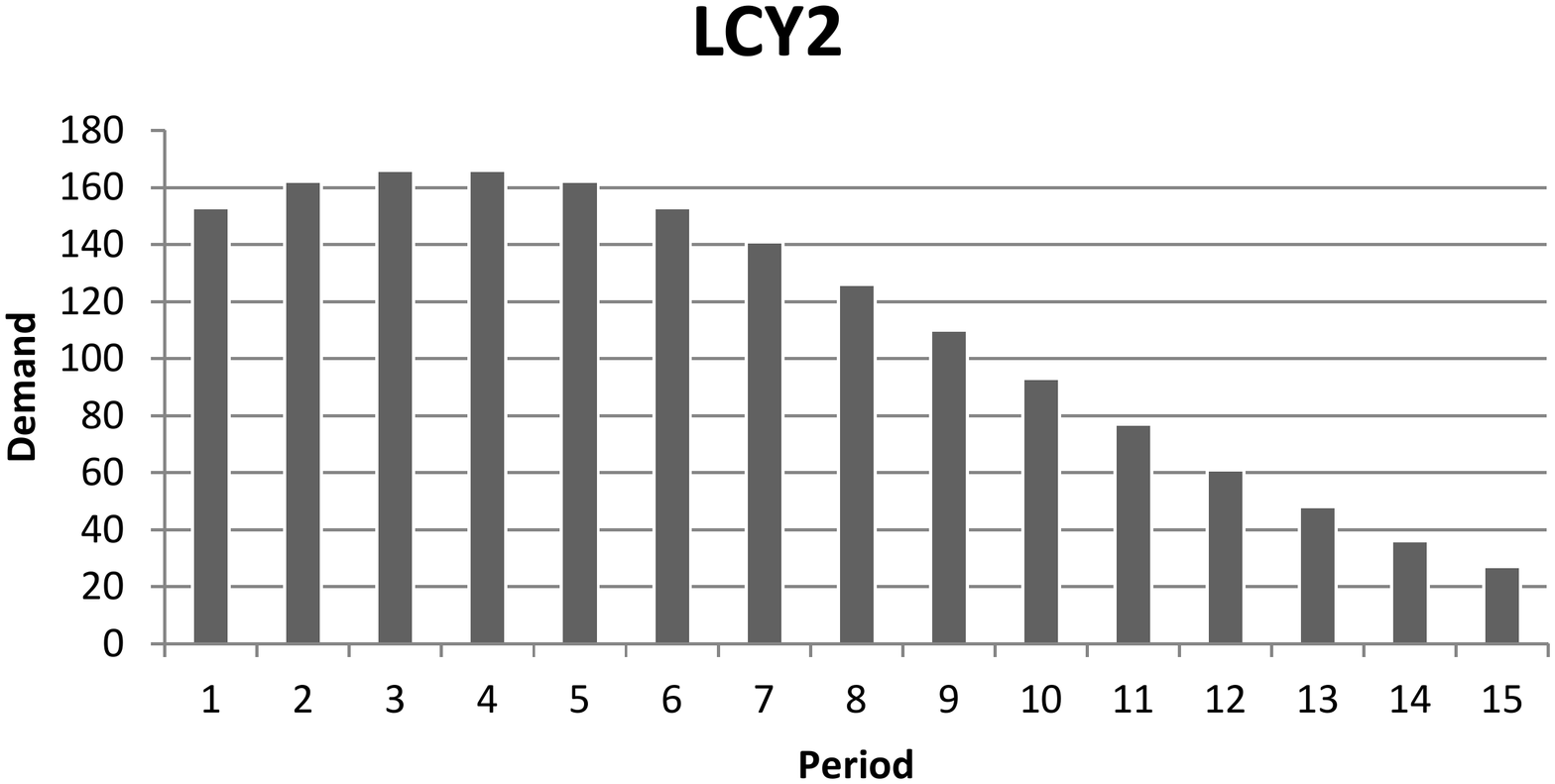}
                
        \end{subfigure}
        \begin{subfigure}[b]{0.45\textwidth}
                \includegraphics[width=\linewidth]{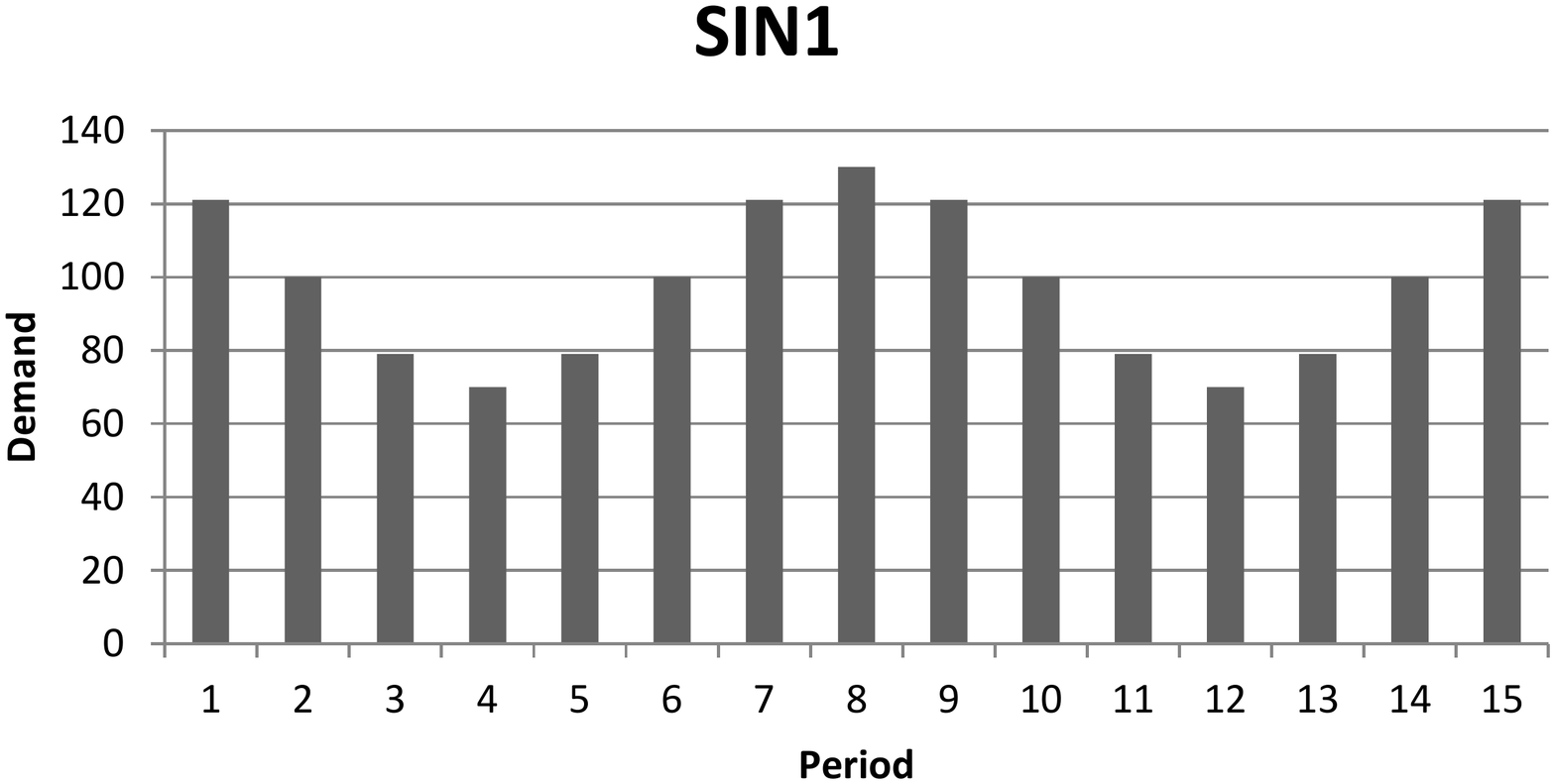}
               
        \end{subfigure}
        \begin{subfigure}[b]{0.45\textwidth}
                \includegraphics[width=\linewidth]{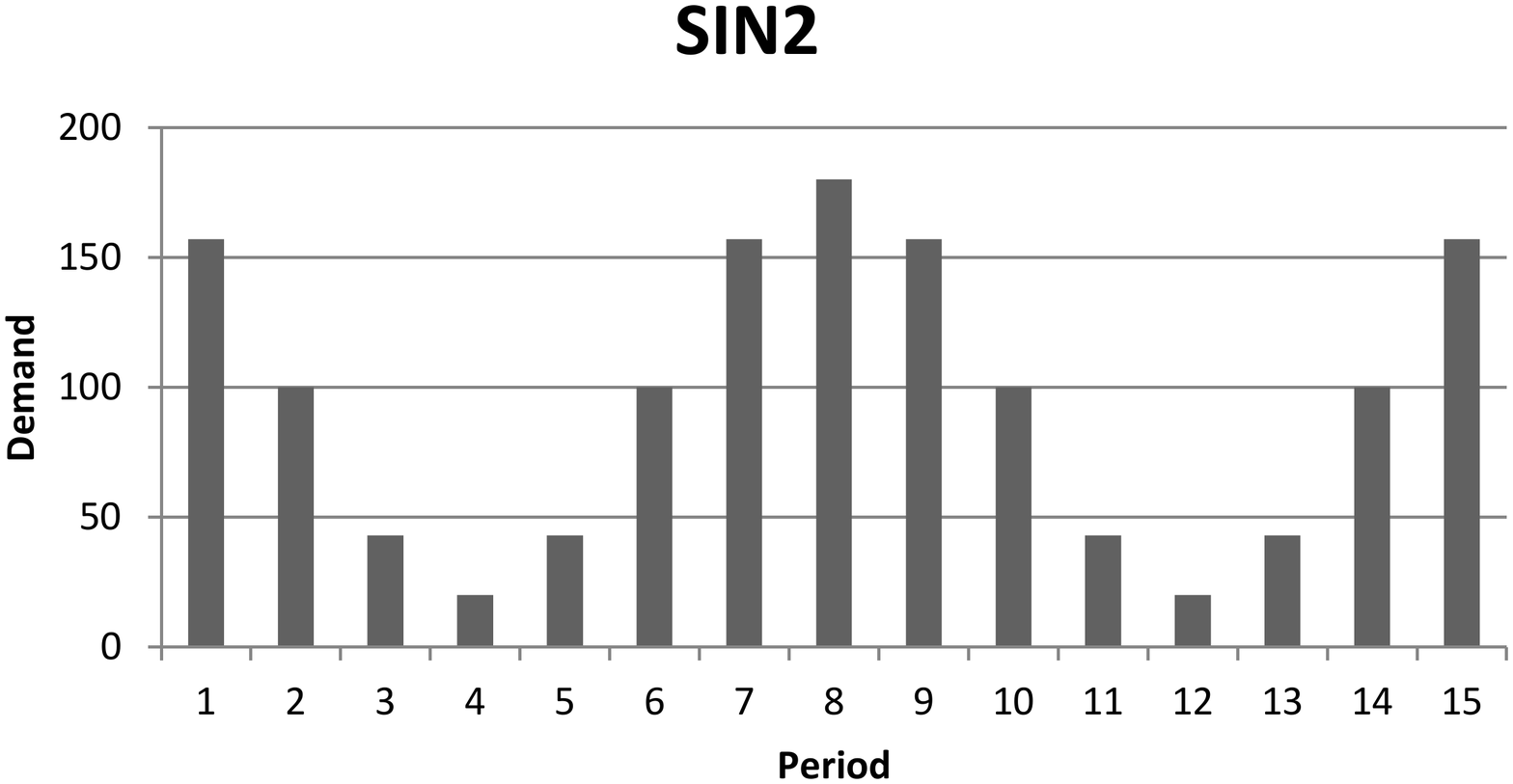}
                
        \end{subfigure}
        \begin{subfigure}[b]{0.45\textwidth}
                \includegraphics[width=\linewidth]{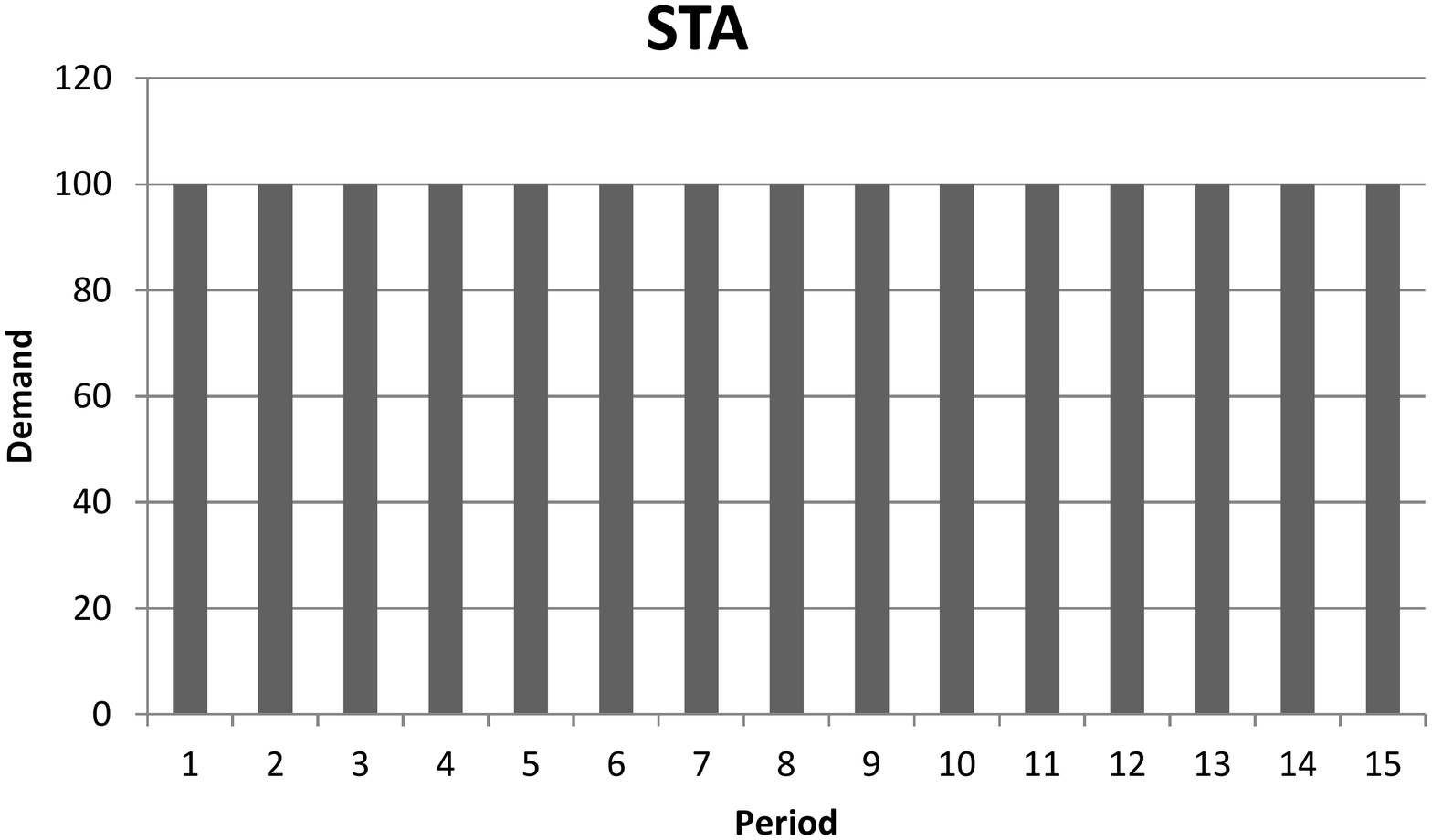}
                
        \end{subfigure}
        \begin{subfigure}[b]{0.45\textwidth}
                \includegraphics[width=\linewidth]{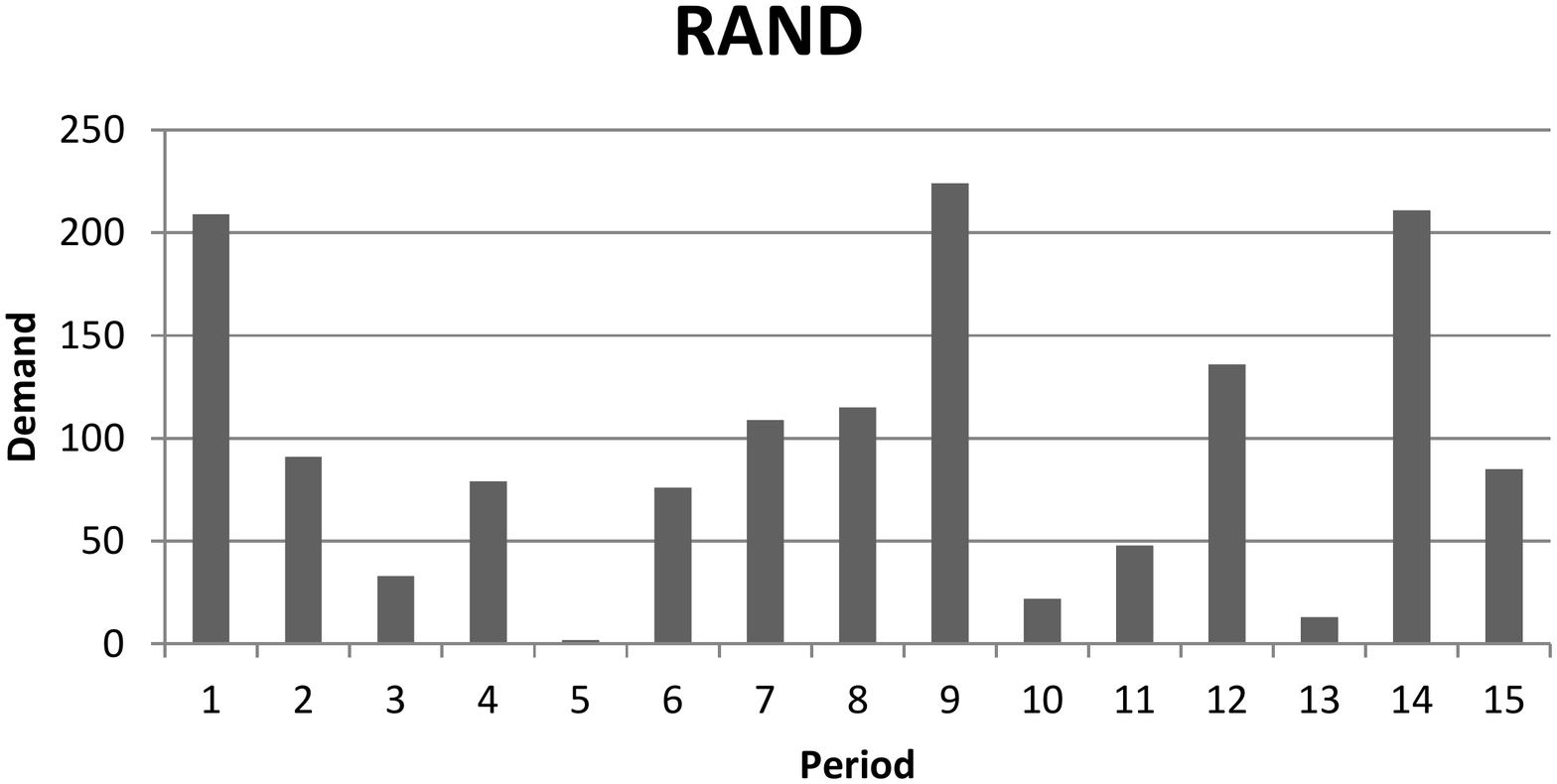}
                
        \end{subfigure}
       
        \caption{Demand Patterns}
        \label{fig:graph}
  
\end{figure*}
\end{footnotesize}

\begin{landscape}
\begin{table}[htbp]

\begin{minipage}{.5\textwidth}
    \centering
    
\begin{footnotesize}

\caption{The combinations of CMILP and exact algorithms of TSP }
\begin{tabular}{cccccccccccccccc}
\hline
\multicolumn{1}{l}{} & \multicolumn{ 3}{c}{CMILP+MTZ} & \multicolumn{ 3}{c}{CMILP+MTZ+2Clq} & \multicolumn{ 3}{c}{CMILP+DL} & \multicolumn{ 3}{c}{CMILP+SC} & \multicolumn{ 3}{c}{CMILP+2C } \\ 
\hline
 & Time (sec) & \%Gap & Best solution & Time (sec) & \%Gap & Best solution & Time (sec) & \%Gap & Best solution & Time (sec) & \%Gap & Best solution & Time (sec) & \%Gap & Best solution \\ 
 \hline
Scenario 1 & 3600.00 & 42.60 & 31000 & 3600.00 & 40.16 & 29740 & 3600.00 & 40.63 & 30360 & 3600.00 & 0.33 & 28740 & 3600.00 & 13.85 & 29120 \\ 
Scenario 2 & 3600.00 & 41.83 & 28435 & 3600.00 & 38.38 & 27035 & 3600.00 & 46.15 & 30115 & 190.13 & 0.00 & 26265 & 3600.00 & 13.12 & 26435 \\ 
Scenario 3 & 3600.00 & 42.55 & 30491 & 3600.00 & 39.64 & 29331 & 3600.00 & 38.85 & 29120 & 241.54 & 0.00 & 28027 & 3600.00 & 11.54 & 28145 \\ 
Scenario 4 & 3600.00 & 1.26 & 103020 & 3600.00 & 1.69 & 103020 & 3600.00 & 1.37 & 103020 & 132.01 & 0.00 & 103020 & 3600.00 & 0.74 & 103020 \\ 
Scenario 5 & 3600.00 & 2.30 & 90220 & 3600.00 & 2.70 & 90785 & 3600.00 & 1.87 & 90220 & 149.13 & 0.00 & 90220 & 3600.00 & 0.95 & 90220 \\ 
Scenario 6 & 3600.00 & 1.95 & 102119 & 3600.00 & 1.94 & 102119 & 3600.00 & 1.86 & 102119 & 186.08 & 0.00 & 102119 & 3600.00 & 1.10 & 102328 \\ 
Scenario 7 & 3600.00 & 4.16 & 111780 & 3600.00 & 3.06 & 110860 & 3600.00 & 3.52 & 111420 & 245.04 & 0.00 & 110520 & 3600.00 & 1.81 & 111140 \\ 
Scenario 8 & 3600.00 & 2.46 & 94775 & 3600.00 & 2.86 & 94975 & 3600.00 & 2.67 & 94975 & 199.46 & 0.00 & 94775 & 3600.00 & 1.04 & 94775 \\ 
Scenario 9 & 3600.00 & 3.05 & 108239 & 3600.00 & 2.92 & 108132 & 3600.00 & 2.85 & 108019 & 218.36 & 0.00 & 107893 & 3600.00 & 1.72 & 107989 \\ 
Scenario 10 & 3600.00 & 81.92 & 101860 & 3600.00 & 81.33 & 101620 & 3600.00 & 83.36 & 113080 & 168.49 & 0.00 & 100020 & 3600.00 & 0.65 & 100020 \\ 
Scenario 11 & 3600.00 & 81.91 & 92390 & 3600.00 & 82.62 & 98025 & 3600.00 & 80.89 & 88625 & 154.35 & 0.00 & 87270 & 3600.00 & 0.46 & 87270 \\ 
Scenario 12 & 3600.00 & 83.18 & 110143 & 3600.00 & 82.09 & 104604 & 3600.00 & 84.7 & 118467 & 423.20 & 0.00 & 99808 & 3600.00 & 1.21 & 99808 \\ 
Scenario 13 & 3600.00 & 31.23 & 148020 & 3600.00 & 31.48 & 148020 & 3600.00 & 31.19 & 148020 & 134.97 & 0.00 & 145080 & 3600.00 & 0.16 & 145080 \\ 
Scenario 14 & 3600.00 & 32.82 & 131305 & 3600.00 & 30.37 & 126605 & 3600.00 & 32.94 & 131305 & 87.31 & 0.00 & 124255 & 3600.00 & 0.10 & 124255 \\ 
Scenario 15 & 3600.00 & 31.61 & 145898 & 3600.00 & 34.72 & 153045 & 3600.00 & 30.89 & 144459 & 199.41 & 0.00 & 144279 & 3600.00 & 0.10 & 144279 \\ 
Scenario 16 & 3600.00 & 31.47 & 157180 & 3600.00 & 32.94 & 160440 & 3600.00 & 28.28 & 150080 & 140.67 & 0.00 & 150080 & 3600.00 & 0.25 & 150080 \\ 
Scenario 17 & 3600.00 & 36.18 & 143355 & 3600.00 & 36.04 & 143355 & 3600.00 & 36.08 & 143355 & 116.10 & 0.00 & 129255 & 1970.35 & 0.00 & 129255 \\ 
Scenario 18 & 3600.00 & 31.30 & 152744 & 3600.00 & 31.44 & 152804 & 3600.00 & 30.43 & 150838 & 161.46 & 0.00 & 149279 & 3600.00 & 0.15 & 149279 \\ 
\hline
\end{tabular}
\label{table:CMILP}
\end{footnotesize}

\end{minipage}

\vspace*{1.5 cm}
\begin{minipage}{.5\textwidth}
    \centering
    
    \begin{footnotesize}
    
\caption{The combinations of SP and exact exact algorithms of TSP }
\begin{tabular}{cccccccccccccccc}
\hline
 \multicolumn{1}{l}{} & \multicolumn{ 3}{c}{SP+MTZ} & \multicolumn{ 3}{c}{SP+MTZ+2Clq} & \multicolumn{ 3}{c}{SP+DL} & \multicolumn{ 3}{c}{SP+SC} & \multicolumn{ 3}{c}{SP+2C } \\ 
 \hline
 & Time (sec) & \%Gap & Best solution & Time (sec) & \%Gap & Best solution & Time (sec) & \%Gap & Best solution & Time (sec) & \%Gap & Best solution & Time (sec) & \%Gap & Best solution \\ 
 \hline
Scenario 1 & 3600.00 & 55.15 & 34400.00 & 3600.00 & 48.54 & 30140.00 & 3600.00 & 44.69 & 30380.00 & 3600.00 & 30.83 & 30220.00 & 3600.00 & 26.77 & 30040.00 \\ 
Scenario 2 & 3600.00 & 46.64 & 28270.00 & 3600.00 & 44.82 & 27755.00 & 3600.00 & 47.27 & 29700.00 & 3600.00 & 27.59 & 26940.00 & 3600.00 & 28.86 & 28055.00 \\ 
Scenario 3 & 3600.00 & 54.20 & 31808.00 & 3600.00 & 47.48 & 30349.00 & 3600.00 & 44.25 & 29789.00 & 3600.00 & 26.69 & 29187.00 & 3600.00 & 30.03 & 29130.00 \\ 
Scenario 4 & 3600.00 & 1.87 & 103020.00 & 3600.00 & 1.86 & 103020.00 & 3600.00 & 1.69 & 103020.00 & 3600.00 & 0.48 & 103020.00 & 3600.00 & 0.94 & 103020.00 \\ 
Scenario 5 & 3600.00 & 2.23 & 90220.00 & 3600.00 & 2.44 & 90220.00 & 3600.00 & 2.55 & 90220.00 & 3600.00 & 0.47 & 90220.00 & 3600.00 & 1.00 & 90220.00 \\ 
Scenario 6 & 3600.00 & 2.51 & 102179.00 & 3600.00 & 2.55 & 102119.00 & 3600.00 & 2.42 & 102119.00 & 3600.00 & 0.75 & 102119.00 & 3600.00 & 1.72 & 102119.00 \\ 
Scenario 7 & 3600.00 & 3.40 & 110860.00 & 3600.00 & 3.14 & 110760.00 & 3600.00 & 3.17 & 110980.00 & 3600.00 & 0.90 & 110520.00 & 3600.00 & 2.05 & 110860.00 \\ 
Scenario 8 & 3600.00 & 2.53 & 94775.00 & 3600.00 & 2.41 & 94775.00 & 3600.00 & 2.83 & 94775.00 & 3600.00 & 0.49 & 94775.00 & 3600.00 & 1.60 & 94775.00 \\ 
Scenario 9 & 3600.00 & 3.41 & 107989.00 & 3600.00 & 4.03 & 108385.00 & 3600.00 & 4.28 & 108955.00 & 3600.00 & 1.28 & 107893.00 & 3600.00 & 1.67 & 107893.00 \\ 
Scenario 10 & 3600.00 & 86.91 & 127700.00 & 3600.00 & 88.04 & 136980.00 & 3600.00 & 87.48 & 141600.00 & 3600.00 & 0.59 & 100020.00 & 3600.00 & 0.76 & 100020.00 \\ 
Scenario 11 & 3600.00 & 83.01 & 93515.00 & 3600.00 & 84.36 & 101725.00 & 3600.00 & 83.70 & 99195.00 & 3600.00 & 0.29 & 87270.00 & 3600.00 & 1.13 & 87270.00 \\ 
Scenario 12 & 3600.00 & 85.87 & 117240.00 & 3600.00 & 88.22 & 136899.00 & 3600.00 & 85.51 & 119571.00 & 3600.00 & 1.83 & 99808.00 & 3600.00 & 0.64 & 99808.00 \\ 
Scenario 13 & 3600.00 & 31.20 & 146520.00 & 3600.00 & 34.63 & 154680.00 & 3600.00 & 31.97 & 148200.00 & 2875.52 & 0.00 & 145080.00 & 3600.00 & 0.28 & 145080.00 \\ 
Scenario 14 & 3600.00 & 33.11 & 131365.00 & 3600.00 & 33.02 & 131485.00 & 3600.00 & 30.78 & 127085.00 & 111.88 & 0.00 & 124255.00 & 3600.00 & 0.31 & 124255.00 \\ 
Scenario 15 & 3600.00 & 35.50 & 153853.00 & 3600.00 & 32.75 & 147938.00 & 3600.00 & 32.22 & 146438.00 & 1615.71 & 0.00 & 144279.00 & 3600.00 & 0.04 & 144279.00 \\ 
Scenario 16 & 3600.00 & 29.52 & 152600.00 & 3600.00 & 31.87 & 157900.00 & 3600.00 & 31.55 & 156740.00 & 879.27 & 0.00 & 150080.00 & 3600.00 & 0.40 & 150080.00 \\ 
Scenario 17 & 3600.00 & 29.36 & 130035.00 & 3600.00 & 29.47 & 130155.00 & 3600.00 & 29.24 & 129735.00 & 1635.82 & 0.00 & 129255.00 & 3600.00 & 0.00 & 129255.00 \\ 
Scenario 18 & 3600.00 & 30.82 & 151138.00 & 3600.00 & 31.42 & 152398.00 & 3600.00 & 30.90 & 151378.00 & 3600.00 & 0.04 & 149279.00 & 3438.31 & 0.00 & 149279.00 \\ 
\hline
\end{tabular}
\label{table:SP}
\end{footnotesize}

\end{minipage}
\end{table}
\end{landscape}

\begin{table*}[htbp]
\center
\begin{footnotesize}

\caption{LP relaxation of CMILP+TSP combinations}
\begin{tabular}{lcccccccccc}
\hline
 & \multicolumn{ 2}{c}{CMILP+MTZ} & \multicolumn{ 2}{c}{CMILP+MTZ+2Clq} & \multicolumn{ 2}{c}{CMILP+DL} & \multicolumn{ 2}{c}{CMILP+SC} & \multicolumn{ 2}{c}{CMILP+2C } \\ 
 \hline
 & \multicolumn{1}{l}{Time (sec)} & \multicolumn{1}{l}{\%Gap} & \multicolumn{1}{l}{Time (sec)} & \multicolumn{1}{l}{\%Gap} & \multicolumn{1}{l}{Time (sec)} & \multicolumn{1}{l}{\%Gap} & \multicolumn{1}{l}{Time (sec)} & \multicolumn{1}{l}{\%Gap} & \multicolumn{1}{l}{Time (sec)} & \multicolumn{1}{l}{\%Gap} \\ 
 \hline
Scenario 1 & 28.50 & 0.93 & 32.82 & 0.62 & 30.47 & 0.62 & 29.80 & 0.62 & 31.98 & 0.62 \\ 
Scenario 2 & 28.78 & 0.93 & 29.39 & 0.59 & 28.38 & 0.59 & 30.45 & 0.59 & 33.20 & 0.59 \\ 
Scenario 3 & 27.98 & 0.93 & 30.29 & 0.61 & 29.26 & 0.61 & 30.61 & 0.64 & 34.00 & 0.64 \\ 
Scenario 4 & 34.07 & 0.06 & 42.37 & 0.04 & 40.41 & 0.04 & 54.90 & 0.05 & 74.40 & 0.05 \\ 
Scenario 5 & 46.83 & 0.07 & 47.45 & 0.04 & 51.04 & 0.04 & 93.50 & 0.05 & 66.15 & 0.05 \\ 
Scenario 6 & 72.37 & 0.08 & 49.84 & 0.04 & 82.76 & 0.04 & 222.77 & 0.06 & 110.68 & 0.06 \\ 
Scenario 7 & 60.86 & 0.07 & 49.15 & 0.04 & 54.49 & 0.04 & 282.69 & 0.05 & 151.99 & 0.05 \\ 
Scenario 8 & 52.07 & 0.08 & 49.23 & 0.05 & 50.50 & 0.05 & 119.69 & 0.06 & 67.32 & 0.06 \\ 
Scenario 9 & 87.54 & 0.08 & 59.45 & 0.05 & 70.93 & 0.05 & 251.33 & 0.06 & 214.50 & 0.06 \\ 
Scenario 10 & 28.95 & 0.95 & 31.53 & 0.87 & 29.13 & 0.87 & 34.17 & 0.05 & 45.02 & 0.05 \\ 
Scenario 11 & 28.95 & 0.95 & 31.14 & 0.85 & 30.48 & 0.85 & 32.46 & 0.06 & 36.59 & 0.06 \\ 
Scenario 12 & 29.21 & 0.95 & 31.08 & 0.87 & 30.93 & 0.87 & 37.73 & 0.05 & 35.18 & 0.05 \\ 
Scenario 13 & 34.69 & 0.33 & 53.46 & 0.31 & 62.54 & 0.31 & 37.50 & 0.02 & 39.96 & 0.02 \\ 
Scenario 14 & 42.15 & 0.32 & 41.15 & 0.30 & 50.80 & 0.30 & 32.85 & 0.03 & 33.30 & 0.03 \\ 
Scenario 15 & 69.61 & 0.34 & 43.10 & 0.32 & 78.40 & 0.32 & 36.80 & 0.02 & 36.36 & 0.02 \\ 
Scenario 16 & 96.88 & 0.31 & 49.32 & 0.29 & 54.45 & 0.29 & 42.30 & 0.02 & 35.11 & 0.02 \\ 
Scenario 17 & 53.57 & 0.32 & 45.97 & 0.30 & 58.71 & 0.30 & 37.59 & 0.02 & 33.40 & 0.02 \\ 
Scenario 18 & 119.05 & 0.33 & 48.88 & 0.31 & 56.04 & 0.31 & 41.63 & 0.02 & 35.51 & 0.02 \\ 
\hline
\end{tabular}
\label{table:LPCMILP}
\end{footnotesize}
\end{table*}

\begin{table*}[htbp]
\center
\begin{footnotesize}

\caption{LP relaxation of SP+TSP combinations}
\begin{tabular}{lcccccccccc}
\hline
& \multicolumn{ 2}{c}{SP+MTZ} & \multicolumn{ 2}{c}{SP+MTZ+2Clq} & \multicolumn{ 2}{c}{SP+DL} & \multicolumn{ 2}{c}{SP+SC} & \multicolumn{ 2}{c}{SP+2C } \\ 
\hline
 & Time (sec) & \%Gap & Time (sec) & \%Gap & Time (sec) & \%Gap & Time (sec) & \%Gap & Time (sec) & \%Gap \\ 
 \hline
Scenario 1 & 39.94 & 0.98 & 42.80 & 0.86 & 39.89 & 0.86 & 40.18 & 0.62 & 43.87 & 0.62 \\ 
Scenario 2 & 41.27 & 0.98 & 42.13 & 0.86 & 41.08 & 0.86 & 39.39 & 0.59 & 44.37 & 0.59 \\ 
Scenario 3 & 39.93 & 0.98 & 42.03 & 0.85 & 41.21 & 0.85 & 40.13 & 0.64 & 45.10 & 0.64 \\ 
Scenario 4 & 40.82 & 0.07 & 44.06 & 0.06 & 42.62 & 0.06 & 43.38 & 0.05 & 47.01 & 0.05 \\ 
Scenario 5 & 42.93 & 0.08 & 46.39 & 0.06 & 42.35 & 0.06 & 40.55 & 0.06 & 46.65 & 0.06 \\ 
Scenario 6 & 42.53 & 0.09 & 45.83 & 0.08 & 43.40 & 0.08 & 41.31 & 0.07 & 47.30 & 0.07 \\ 
Scenario 7 & 42.49 & 0.08 & 48.25 & 0.07 & 43.40 & 0.07 & 43.86 & 0.06 & 46.53 & 0.06 \\ 
Scenario 8 & 42.37 & 0.09 & 43.55 & 0.08 & 43.47 & 0.08 & 39.99 & 0.07 & 46.58 & 0.07 \\ 
Scenario 9 & 42.62 & 0.09 & 43.86 & 0.08 & 43.50 & 0.08 & 39.02 & 0.07 & 46.64 & 0.07 \\ 
Scenario 10 & 41.95 & 0.96 & 40.60 & 0.93 & 43.49 & 0.93 & 39.16 & 0.05 & 44.58 & 0.05 \\ 
Scenario 11 & 40.49 & 0.97 & 42.10 & 0.93 & 41.46 & 0.93 & 39.37 & 0.06 & 45.03 & 0.06 \\ 
Scenario 12 & 40.54 & 0.97 & 43.30 & 0.93 & 39.80 & 0.93 & 38.38 & 0.08 & 45.72 & 0.08 \\ 
Scenario 13 & 42.67 & 0.33 & 44.11 & 0.32 & 42.59 & 0.32 & 40.02 & 0.02 & 46.80 & 0.02 \\ 
Scenario 14 & 42.62 & 0.32 & 43.33 & 0.32 & 41.61 & 0.32 & 40.34 & 0.03 & 46.82 & 0.03 \\ 
Scenario 15 & 43.06 & 0.35 & 44.43 & 0.34 & 43.71 & 0.34 & 40.53 & 0.04 & 47.26 & 0.04 \\ 
Scenario 16 & 42.41 & 0.31 & 43.51 & 0.31 & 41.71 & 0.31 & 39.97 & 0.03 & 47.08 & 0.03 \\ 
Scenario 17 & 42.44 & 0.33 & 43.96 & 0.32 & 43.17 & 0.32 & 40.14 & 0.03 & 46.36 & 0.03 \\ 
Scenario 18 & 42.70 & 0.34 & 44.49 & 0.33 & 41.88 & 0.33 & 41.51 & 0.05 & 46.62 & 0.05 \\ 
\hline
\end{tabular}
\label{table:LPSP}
\end{footnotesize}
\end{table*}

\begin{table*}[htbp]
\center
\begin{footnotesize}

\caption{The Optimality gap(\%) for the combinations of CMILP+TSP algorithms proposed by \cite{bektacs2014}}
\begin{tabular}{lccccccccc}
\hline
\multicolumn{ 1}{c}{} & CMILP+ & CMILP+ & CMILP+ & CMILP+ & CMILP+ & CMILP+ & CMILP+ & CMILP+ & CMILP+DL+ \\ 
\multicolumn{ 1}{c}{} & DL+3CLQ & DL+NR &  DL+L3 & DL+2P & DL+R & R+2P & NR+2P & NR+R+2P &  NR+R+2P \\ 
\hline
Scenario 1 & 0.38 & 0.38 & 0.45 & 0.44 & 0.38 & 0.38 & 0.38 & 0.38 & 0.38 \\ 
Scenario 2 & 0.38 & 0.38 & 0.44 & 0.43 & 0.38 & 0.38 & 0.38 & 0.38 & 0.38 \\ 
Scenario 3 & 0.37 & 0.37 & 0.44 & 0.43 & 0.37 & 0.37 & 0.37 & 0.37 & 0.37 \\ 
Scenario 4 & 0.02 & 0.02 & 0.03 & 0.03 & 0.02 & 0.02 & 0.02 & 0.02 & 0.02 \\ 
Scenario 5 & 0.03 & 0.03 & 0.03 & 0.03 & 0.03 & 0.03 & 0.03 & 0.03 & 0.03 \\ 
Scenario 6 & 0.02 & 0.02 & 0.03 & 0.02 & 0.02 & 0.02 & 0.02 & 0.02 & 0.02 \\ 
Scenario 7 & 0.03 & 0.03 & 0.03 & 0.03 & 0.03 & 0.03 & 0.03 & 0.03 & 0.03 \\ 
Scenario 8 & 0.03 & 0.03 & 0.04 & 0.04 & 0.03 & 0.03 & 0.03 & 0.03 & 0.03 \\ 
Scenario 9 & 0.03 & 0.03 & 0.03 & 0.03 & 0.03 & 0.03 & 0.03 & 0.03 & 0.03 \\ 
Scenario 10 & 0.82 & 0.82 & 0.83 & 0.83 & 0.82 & 0.82 & 0.82 & 0.82 & 0.82 \\ 
Scenario 11 & 0.81 & 0.81 & 0.82 & 0.82 & 0.81 & 0.81 & 0.81 & 0.81 & 0.81 \\ 
Scenario 12 & 0.82 & 0.82 & 0.83 & 0.83 & 0.82 & 0.82 & 0.82 & 0.82 & 0.82 \\ 
Scenario 13 & 0.30 & 0.30 & 0.31 & 0.31 & 0.30 & 0.30 & 0.30 & 0.30 & 0.30 \\ 
Scenario 14 & 0.29 & 0.29 & 0.29 & 0.29 & 0.29 & 0.29 & 0.29 & 0.29 & 0.29 \\ 
Scenario 15 & 0.30 & 0.30 & 0.31 & 0.30 & 0.30 & 0.30 & 0.30 & 0.30 & 0.30 \\ 
Scenario 16 & 0.28 & 0.28 & 0.28 & 0.28 & 0.28 & 0.28 & 0.28 & 0.28 & 0.28 \\ 
Scenario 17 & 0.29 & 0.29 & 0.29 & 0.29 & 0.29 & 0.29 & 0.29 & 0.29 & 0.29 \\ 
Scenario 18 & 0.29 & 0.29 & 0.30 & 0.29 & 0.29 & 0.29 & 0.29 & 0.29 & 0.29 \\ 
\hline
\end{tabular}
\label{table:gapcmılp}
\end{footnotesize}
\end{table*}

\begin{table*}[htbp]
\center
\begin{footnotesize}
\caption{The Computational times (in seconds) for the combinations of CMILP+TSP algorithms proposed by \cite{bektacs2014}}
\begin{tabular}{lccccccccc}
\hline
\multicolumn{ 1}{c}{} & CMILP+ & CMILP+ & CMILP+ & CMILP+ & CMILP+ & CMILP+ & CMILP+ & CMILP+ & CMILP+DL+ \\ 
\multicolumn{ 1}{c}{} & DL+3CLQ & DL+NR &  DL+L3 & DL+2P & DL+R & R+2P & NR+2P & NR+R+2P &  NR+R+2P \\ 
\hline
Scenario 1 & 75.45 & 78.29 & 99.69 & 1843.28 & 94.02 & 1068.07 & 1007.02 & 1102.41 & 1179.68 \\ 
Scenario 2 & 84.46 & 93.12 & 97.34 & 1015.32 & 92.38 & 785.65 & 732.02 & 806.16 & 1106.99 \\ 
Scenario 3 & 74.43 & 80.65 & 99.15 & 945.24 & 84.08 & 1314.85 & 851.35 & 1089.17 & 841.18 \\ 
Scenario 4 & 83.80 & 117.16 & 99.31 & 404.86 & 100.57 & 583.09 & 432.10 & 553.32 & 549.28 \\ 
Scenario 5 & 92.67 & 103.63 & 206.72 & 679.27 & 116.69 & 578.77 & 491.48 & 661.66 & 530.19 \\ 
Scenario 6 & 140.14 & 291.42 & 1656.04 & 3600.00 & 171.18 & 1556.72 & 1972.35 & 1532.07 & 2355.09 \\ 
Scenario 7 & 105.10 & 122.92 & 382.10 & 962.92 & 108.51 & 716.02 & 533.80 & 733.65 & 638.82 \\ 
Scenario 8 & 99.67 & 109.36 & 231.97 & 1443.81 & 96.63 & 644.46 & 473.48 & 667.46 & 751.39 \\ 
Scenario 9 & 168.84 & 213.07 & 2531.54 & 3600.00 & 204.64 & 1610.78 & 1894.70 & 2014.20 & 1226.56 \\ 
Scenario 10 & 70.23 & 79.24 & 102.26 & 3600.00 & 79.84 & 668.84 & 747.70 & 657.56 & 702.15 \\ 
Scenario 11 & 78.94 & 110.33 & 125.01 & 1289.42 & 96.20 & 808.04 & 728.27 & 711.57 & 909.68 \\ 
Scenario 12 & 75.30 & 84.28 & 126.38 & 2600.50 & 85.05 & 584.08 & 581.21 & 631.57 & 659.13 \\ 
Scenario 13 & 98.09 & 110.87 & 139.36 & 403.24 & 87.18 & 464.86 & 480.90 & 545.68 & 575.52 \\ 
Scenario 14 & 90.75 & 111.42 & 215.81 & 584.06 & 95.43 & 648.40 & 498.00 & 581.33 & 573.13 \\ 
Scenario 15 & 426.08 & 327.40 & 1646.94 & 3600.00 & 457.77 & 2308.55 & 3205.38 & 3600.00 & 2722.97 \\ 
Scenario 16 & 88.15 & 98.21 & 941.78 & 3036.53 & 103.56 & 600.14 & 715.69 & 726.25 & 755.38 \\ 
Scenario 17 & 89.60 & 122.05 & 295.00 & 1875.84 & 107.27 & 612.66 & 717.07 & 578.27 & 725.89 \\ 
Scenario 18 & 348.92 & 324.75 & 3600.00 & 3600.00 & 397.32 & 3600.00 & 3600.00 & 3600.00 & 3600.00 \\ 
\hline
\end{tabular}
\label{table:CMILPbektas}
\end{footnotesize}
\end{table*}

\begin{table*}[htbp]
\center
\begin{footnotesize}

\caption{The optimality gap(\%) for the combinations of SP+TSP algorithms proposed by \cite{bektacs2014}}
\begin{tabular}{lccccccccc}
\hline
\multicolumn{1}{c}{} & \multicolumn{ 1}{c}{SP+} & SP+ & SP+ & SP+ & SP+ & SP+ & SP+ & SP+ & SP+DL+ \\ 
\multicolumn{1}{c}{} & \multicolumn{ 1}{c}{DL+3CLQ   } & DL+NR &  DL+L3 & DL+2P & DL+R & R+2P & NR+2P & NR+R+2P &  NR+R+2P \\ 
\hline
Scenario 1 & 0.80 & 0.80 & 0.82 & 0.82 & 0.80 & 0.80 & 0.80 & 0.80 & 0.80 \\ 
Scenario 2 & 0.81 & 0.81 & 0.83 & 0.82 & 0.81 & 0.81 & 0.81 & 0.81 & 0.81 \\ 
Scenario 3 & 0.79 & 0.79 & 0.81 & 0.81 & 0.79 & 0.79 & 0.79 & 0.79 & 0.79 \\ 
Scenario 4 & 0.05 & 0.05 & 0.05 & 0.05 & 0.05 & 0.05 & 0.05 & 0.05 & 0.05 \\ 
Scenario 5 & 0.06 & 0.06 & 0.06 & 0.06 & 0.06 & 0.06 & 0.06 & 0.06 & 0.06 \\ 
Scenario 6 & 0.07 & 0.07 & 0.07 & 0.07 & 0.07 & 0.07 & 0.07 & 0.07 & 0.07 \\ 
Scenario 7 & 0.06 & 0.06 & 0.07 & 0.07 & 0.06 & 0.06 & 0.06 & 0.06 & 0.06 \\ 
Scenario 8 & 0.07 & 0.07 & 0.07 & 0.07 & 0.07 & 0.07 & 0.07 & 0.07 & 0.07 \\ 
Scenario 9 & 0.07 & 0.07 & 0.08 & 0.08 & 0.07 & 0.07 & 0.07 & 0.07 & 0.07 \\ 
Scenario 10 & 0.91 & 0.91 & 0.92 & 0.92 & 0.91 & 0.91 & 0.91 & 0.91 & 0.91 \\ 
Scenario 11 & 0.92 & 0.92 & 0.92 & 0.92 & 0.92 & 0.92 & 0.92 & 0.92 & 0.92 \\ 
Scenario 12 & 0.92 & 0.92 & 0.92 & 0.92 & 0.92 & 0.92 & 0.92 & 0.92 & 0.92 \\ 
Scenario 13 & 0.32 & 0.32 & 0.32 & 0.32 & 0.32 & 0.32 & 0.32 & 0.32 & 0.32 \\ 
Scenario 14 & 0.31 & 0.31 & 0.31 & 0.31 & 0.31 & 0.31 & 0.31 & 0.31 & 0.31 \\ 
Scenario 15 & 0.34 & 0.34 & 0.34 & 0.34 & 0.34 & 0.34 & 0.34 & 0.34 & 0.34 \\ 
Scenario 16 & 0.31 & 0.31 & 0.31 & 0.31 & 0.31 & 0.31 & 0.31 & 0.31 & 0.31 \\ 
Scenario 17 & 0.32 & 0.32 & 0.32 & 0.32 & 0.32 & 0.32 & 0.32 & 0.32 & 0.32 \\ 
Scenario 18 & 0.33 & 0.33 & 0.33 & 0.33 & 0.33 & 0.33 & 0.33 & 0.33 & 0.33 \\ 
\hline
\end{tabular}
\label{table:gapSP}
\end{footnotesize}
\end{table*}
\begin{table*}[htbp]
\center
\begin{footnotesize}

\caption{The Computational times the combinations of SP+TSP algorithms proposed by \cite{bektacs2014}}
\begin{tabular}{lccccccccc}
\hline
\multicolumn{1}{c}{} & \multicolumn{ 1}{c}{SP+} & SP+ & SP+ & SP+ & SP+ & SP+ & SP+ & SP+ & SP+DL+ \\ 
\multicolumn{1}{c}{} & \multicolumn{ 1}{c}{DL+3CLQ} & DL+NR &  DL+L3 & DL+2P & DL+R & R+2P & NR+2P & NR+R+2P &  NR+R+2P \\ 
\hline
Scenario 1 & 99.01 & 100.96 & 63.56 & 101.94 & 91.97 & 163.14 & 162.29 & 232.10 & 232.11 \\ 
Scenario 2 & 73.59 & 89.31 & 55.44 & 94.64 & 86.92 & 163.70 & 161.20 & 233.89 & 233.39 \\ 
Scenario 3 & 76.70 & 88.88 & 56.38 & 95.48 & 87.54 & 162.60 & 162.18 & 234.02 & 233.79 \\ 
Scenario 4 & 75.43 & 90.08 & 57.30 & 95.10 & 89.39 & 164.44 & 162.33 & 237.56 & 235.32 \\ 
Scenario 5 & 75.75 & 90.79 & 57.13 & 96.28 & 89.09 & 165.49 & 162.92 & 241.00 & 234.40 \\ 
Scenario 6 & 75.23 & 90.41 & 57.11 & 99.28 & 89.02 & 167.30 & 165.85 & 236.17 & 234.28 \\ 
Scenario 7 & 76.15 & 91.19 & 58.00 & 96.86 & 88.63 & 162.82 & 162.20 & 232.74 & 240.71 \\ 
Scenario 8 & 76.61 & 89.99 & 57.72 & 96.49 & 88.59 & 163.16 & 160.94 & 235.30 & 238.08 \\ 
Scenario 9 & 77.21 & 94.47 & 57.99 & 97.31 & 91.04 & 163.12 & 161.66 & 233.00 & 237.38 \\ 
Scenario 10 & 74.50 & 94.04 & 56.49 & 94.77 & 88.29 & 162.72 & 159.88 & 297.86 & 302.65 \\ 
Scenario 11 & 74.04 & 88.95 & 55.90 & 94.01 & 87.88 & 161.29 & 159.53 & 301.19 & 296.99 \\ 
Scenario 12 & 74.64 & 100.79 & 56.04 & 94.61 & 87.50 & 160.86 & 160.88 & 311.06 & 303.91 \\ 
Scenario 13 & 75.89 & 88.23 & 57.76 & 97.21 & 88.63 & 163.77 & 161.68 & 299.14 & 307.28 \\ 
Scenario 14 & 74.92 & 88.01 & 57.05 & 96.45 & 88.85 & 164.64 & 163.49 & 305.52 & 306.08 \\ 
Scenario 15 & 80.32 & 88.43 & 57.36 & 97.13 & 89.36 & 163.00 & 164.67 & 310.07 & 300.93 \\ 
Scenario 16 & 76.34 & 88.20 & 57.41 & 97.59 & 88.80 & 164.03 & 163.04 & 307.31 & 244.32 \\ 
Scenario 17 & 76.58 & 91.58 & 57.29 & 96.40 & 89.02 & 163.77 & 169.32 & 300.15 & 244.89 \\ 
Scenario 18 & 76.43 & 88.27 & 57.31 & 97.39 & 89.08 & 163.01 & 165.25 & 299.23 & 248.05 \\ 
\hline
\end{tabular}
\label{table:SPbektas}
\end{footnotesize}
\end{table*}

\begin{table*}[htbp]

\begin{footnotesize}
  \caption{The average computational times (in seconds) of combination of CMILP and exact algorithms of TSP on 80 instances}   

\begin{center}
\begin{tabular}{ccccccc}
\hline
n & r & CMILP+MTZ & CMILP+MTZ+2Clq & CMILP+DL & CMILP+SC & CMILP+2C  \\ 
\hline
3 & 5 & 0.43 & 0.78 & 0.46 & 0.48 & 0.50 \\ 
3 & 10 & 0.71 & 0.78 & 0.72 & 0.88 & 1.10 \\ 
3 & 15 & 2.13 & 1.88 & 1.88 & 2.32 & 2.63 \\ 
3 & 20 & 25.96 & 27.39 & 22.46 &  4.82 & 5.89 \\ 
3 & 25 & 227.72 & 255.11 & 385.28 &  9.27 & 11.65 \\ 
3 & 30 & 37.27 & 46.44 & 46.52 &  17.82 & 21.97 \\ 
3 & 35 & 3495.71 & 1363.08 & 1108.21 & 31.65 & 36.34 \\ 
3 & 40 & 120.75 & 127.70 & 154.93 &  53.32 & 59.40 \\ 
3 & 45 & 98.77 & 83.54 & 70.79 &  86.12 & 101.78 \\ 
3 & 50 & 1673.92 & 816.66 & 1425.62 &  131.76 & 160.49 \\ 
\hline
6 & 5 & 0.52 & 0.52 & 0.49 & 0.62 & 0.67 \\ 
6 & 10 & 0.97 & 1.06 & 0.97 & 1.62 & 2.07 \\ 
6 & 15 & 3.06 & 3.26 & 2.94 & 5.81 & 6.18 \\ 
6 & 20 & 53.23 & 48.17 & 43.15 &  14.95 & 15.38 \\ 
6 & 25 & 519.55 & 762.68 & 972.94 &  32.48 & 36.12 \\ 
6 & 30 & 77.10 & 72.25 & 74.25 &  66.70 & 72.07 \\ 
\hline
\end{tabular}
\end{center}
\label{table:tab1}
\end{footnotesize}
\end{table*}

\begin{table*}[htbp]

\begin{footnotesize}

  \caption{The average computational times (in seconds) combination of SP and exact algorithms of TSP on 80 instances}

\begin{center}

\begin{tabular}{ccccccc}
\hline
n & r & SP+MTZ & SP+MTZ+2Clq & SP+DL & SP+SC & SP+2C \\
\hline 
3 & 5 & 0.26 & 0.28 & 0.24 & 0.18 & 0.28 \\ 
3 & 10 & 0.74 & 0.75 & 0.75 & 0.59 & 0.86 \\ 
3 & 15 & 2.02 & 2.28 & 2.29 & 1.77 & 2.37 \\ 
3 & 20 & 28.08 & 28.59 & 25.35 & 4.06 & 6.18 \\ 
3 & 25 & 493.88 & 363.27 & 332.72 &  9.04 & 11.61 \\ 
3 & 30 & 65.49 & 66.73 & 64.61 &  17.96 & 21.28 \\ 
3 & 35 & 1341.61 & 1420.89 & 1157.47 &  30.44 & 37.04 \\ 
3 & 40 & 336.16 & 378.16 & 278.05 &  50.86 & 61.89 \\ 
3 & 45 & 166.06 & 195.10 & 175.75 &  86.43 & 106.65 \\ 
3 & 50 & 3501.72 & 1670.40 & 9706.89 &  138.82 & 155.53 \\ 
\hline
6 & 5 & 0.51 & 0.52 & 0.44 & 0.39 & 0.48 \\ 
6 & 10 & 1.74 & 1.83 & 1.76 & 1.57 & 1.75 \\ 
6 & 15 & 5.75 & 5.69 & 5.79 & 5.10 & 5.96 \\ 
6 & 20 & 54.08 & 56.77 & 55.41 &  14.36 & 15.91 \\ 
6 & 25 & 2769.00 & 892.08 & 1220.27 &  31.84 & 35.04 \\ 
6 & 30 & 144.56 & 142.35 & 131.43 &  62.63 & 70.48 \\ 
\hline
\end{tabular}
\end{center}
\label{table:tab2}
\end{footnotesize}
\end{table*}

\begin{table*}[htbp]

\begin{footnotesize}

\caption{The average optimality gap(\%) between LP relaxation solution (\%) and optimal solution of CMILP and TSP combination on 80 instances}

\begin{center}

\begin{tabular}{ccccccc}
\hline
n & r & CMILP+MTZ & CMILP+MTZ+2Clq & CMILP+DL & CMILP+SC & CMILP+2C  \\ 
\hline
3 & 5 & 0.11 & 0.06 & 0.06 & 0.05 & 0.05 \\ 
3 & 10 & 0.17 & 0.05 & 0.05 & 0.08 & 0.08 \\ 
3 & 15 & 0.19 & 0.06 & 0.06 & 0.13 & 0.13 \\ 
3 & 20 & 0.23 & 0.09 & 0.09 & 0.16 & 0.16 \\ 
3 & 25 & 0.20 & 0.11 & 0.11 & 0.13 & 0.13 \\ 
3 & 30 & 0.16 & 0.05 & 0.05 & 0.12 & 0.12 \\ 
3 & 35 & 0.18 & 0.08 & 0.08 & 0.13 & 0.13 \\ 
3 & 40 & 0.19 & 0.01 & 0.06 & 0.15 & 0.15 \\ 
3 & 45 & 0.17 & 0.05 & 0.05 & 0.13 & 0.13 \\ 
3 & 50 & 0.15 & 0.06 & 0.06 & 0.11 & 0.11 \\
\hline 
6 & 5 & 0.10 & 0.05 & 0.05 & 0.05 & 0.05 \\ 
6 & 10 & 0.15 & 0.05 & 0.05 & 0.07 & 0.07 \\ 
6 & 15 & 0.17 & 0.05 & 0.05 & 0.11 & 0.11 \\ 
6 & 20 & 0.20 & 0.08 & 0.08 & 0.14 & 0.14 \\ 
6 & 25 & 0.17 & 0.09 & 0.09 & 0.11 & 0.11 \\ 
6 & 30 & 0.13 & 0.04 & 0.04 & 0.10 & 0.10 \\ 
\hline
\end{tabular}
\end{center}
\label{table:classicrelax}
\end{footnotesize}
\end{table*}

\begin{table*}[htbp]

\begin{footnotesize}

\caption{The average computational times (in seconds) of LP relaxation of CMILP+TSP combination on 80 instances}

\begin{center}

\begin{tabular}{ccccccc}
\hline
n & r & CMILP+MTZ & CMILP+MTZ+2Clq & CMILP+DL & CMILP+SC & CMILP+2C  \\
\hline 
3 & 5 & 0.43 & 0.45 & 0.47 & 0.54 & 0.59 \\ 
3 & 10 & 0.88 & 0.79 & 0.75 & 1.12 & 1.30 \\ 
3 & 15 & 1.27 & 1.37 & 1.36 & 2.36 & 2.72 \\ 
3 & 20 & 1.92 & 2.25 & 2.00 & 4.85 & 5.23 \\ 
3 & 25 & 2.78 & 3.22 & 3.22 & 9.29 & 10.64 \\ 
3 & 30 & 3.99 & 4.30 & 4.09 & 17.54 & 19.11 \\ 
3 & 35 & 4.70 & 4.82 & 5.23 & 30.65 & 32.83 \\ 
3 & 40 & 6.62 & 7.73 & 6.91 & 51.00 & 54.83 \\ 
3 & 45 & 7.32 & 10.06 & 9.06 & 81.03 & 85.52 \\ 
3 & 50 & 9.36 & 12.08 & 11.16 & 124.67 & 133.41 \\ 
\hline
6 & 5 & 0.54 & 0.50 & 0.49 & 0.60 & 0.69 \\ 
6 & 10 & 1.06 & 1.09 & 1.05 & 1.83 & 2.07 \\ 
6 & 15 & 1.90 & 2.05 & 1.91 & 5.39 & 5.93 \\ 
6 & 20 & 3.11 & 3.63 & 3.24 & 14.38 & 14.93 \\ 
6 & 25 & 4.96 & 5.34 & 5.15 & 32.67 & 33.08 \\ 
6 & 30 & 7.52 & 8.46 & 7.70 & 64.91 & 64.08 \\ 
\hline
\end{tabular}
\end{center}
\label{table:classicrelaxtime}
\end{footnotesize}
\end{table*}

\begin{table*}[htbp]

\begin{footnotesize}
 \caption{The average optimality gap(\%) between LP relaxation solution (\%) and optimal solution of combination SP+TSP on 80 instances}

\begin{center}

\begin{tabular}{ccccccc}
\hline
n & r & SP+MTZ & SP+MTZ+2Clq & SP+DL & SP+SC & SP+2C \\ 
\hline
3 & 5 & 0.11 & 0.06 & 0.06 & 0.05 & 0.05 \\ 
3 & 10 & 0.17 & 0.05 & 0.05 & 0.08 & 0.08 \\ 
3 & 15 & 0.19 & 0.06 & 0.06 & 0.13 & 0.13 \\ 
3 & 20 & 0.23 & 0.09 & 0.09 & 0.16 & 0.16 \\ 
3 & 25 & 0.20 & 0.11 & 0.11 & 0.13 & 0.13 \\ 
3 & 30 & 0.16 & 0.05 & 0.05 & 0.12 & 0.12 \\ 
3 & 35 & 0.18 & 0.08 & 0.08 & 0.13 & 0.13 \\ 
3 & 40 & 0.19 & 0.06 & 0.06 & 0.15 & 0.15 \\ 
3 & 45 & 0.17 & 0.05 & 0.05 & 0.13 & 0.13 \\ 
3 & 50 & 0.15 & 0.06 & 0.06 & 0.11 & 0.11 \\
\hline 
6 & 5 & 0.10 & 0.05 & 0.05 & 0.05 & 0.05 \\ 
6 & 10 & 0.15 & 0.05 & 0.05 & 0.07 & 0.07 \\ 
6 & 15 & 0.17 & 0.05 & 0.05 & 0.11 & 0.11 \\ 
6 & 20 & 0.20 & 0.08 & 0.08 & 0.14 & 0.14 \\ 
6 & 25 & 0.17 & 0.09 & 0.09 & 0.11 & 0.11 \\ 
6 & 30 & 0.13 & 0.04 & 0.04 & 0.10 & 0.10 \\
\hline 
\end{tabular}
\end{center}
\label{table:shortestrelax}
\end{footnotesize}
\end{table*}

\begin{table*}[htbp]

\begin{footnotesize}
\caption{The average computational times (in seconds) of LP relaxation of SP+TSP combination on 80 instances}

\begin{center}

\begin{tabular}{ccccccc}
\hline
n & r & SP+MTZ & SP+MTZ+2Clq & SP+DL & SP+SC & SP+2C \\ 
\hline
3 & 5 & 0.13 & 0.13 & 0.12 & 0.13 & 0.17 \\ 
3 & 10 & 0.53 & 0.41 & 0.40 & 0.41 & 0.59 \\ 
3 & 15 & 1.21 & 1.29 & 1.23 & 1.33 & 1.59 \\ 
3 & 20 & 3.10 & 3.21 & 3.09 & 3.20 & 3.65 \\ 
3 & 25 & 6.89 & 6.93 & 6.79 & 6.95 & 8.00 \\ 
3 & 30 & 13.86 & 13.70 & 13.71 & 13.69 & 15.30 \\ 
3 & 35 & 25.60 & 24.92 & 26.88 & 25.18 & 27.69 \\ 
3 & 40 & 44.75 & 43.31 & 42.67 & 43.53 & 48.62 \\ 
3 & 45 & 70.97 & 71.13 & 69.49 & 74.01 & 78.79 \\ 
3 & 50 & 107.14 & 108.67 & 108.66 & 114.25 & 121.47 \\ 
\hline
6 & 5 & 0.28 & 0.40 & 0.33 & 0.30 & 0.49 \\ 
6 & 10 & 1.24 & 1.36 & 1.34 & 1.38 & 1.64 \\ 
6 & 15 & 4.27 & 4.38 & 4.60 & 4.57 & 5.14 \\ 
6 & 20 & 11.94 & 12.65 & 12.52 & 13.13 & 13.78 \\ 
6 & 25 & 27.80 & 28.40 & 29.18 & 29.88 & 31.84 \\ 
6 & 30 & 57.31 & 62.79 & 59.70 & 61.39 & 64.24 \\
\hline 
\end{tabular}
\end{center}
\label{table:shortestrelaxtime}
\end{footnotesize}
\end{table*}

\begin{table*}[htbp]
\center
\begin{footnotesize}
 \caption{The average optimality gap(\%) for the combinations of CMILP+TSP algorithms proposed by \cite{bektacs2014} on 80 instances}
\begin{tabular}{ccccccccccc}
\hline
\multicolumn{1}{l}{} & \multicolumn{1}{l}{} & \multicolumn{1}{l}{} & \multicolumn{1}{l}{} & \multicolumn{1}{l}{} & \multicolumn{1}{l}{} & \multicolumn{1}{l}{} & \multicolumn{1}{l}{} & \multicolumn{1}{l}{} & \multicolumn{1}{l}{} & \multicolumn{1}{l}{} \\ 
\multicolumn{1}{l}{} & \multicolumn{1}{l}{} & CMILP+ & CMILP+ & CMILP+ & CMILP+ & CMILP+ & CMILP+ & CMILP+ & CMILP+ & CMILP+DL+ \\ 
n & r & DL+3CLQ & DL+NR &  DL+L3 & DL+2P & DL+R & R+2P & NR+2P & NR+R+2P &  NR+R+2P \\ 
\hline
3 & 6 & 0.06 & 0.06 & 0.06 & 0.06 & 0.06 & 0.06 & 0.06 & 0.06 & 0.06 \\ 
3 & 11 & 0.05 & 0.05 & 0.05 & 0.05 & 0.05 & 0.05 & 0.05 & 0.05 & 0.05 \\ 
3 & 16 & 0.04 & 0.04 & 0.04 & 0.04 & 0.04 & 0.04 & 0.04 & 0.04 & 0.04 \\ 
3 & 21 & 0.06 & 0.06 & 0.06 & 0.06 & 0.06 & 0.06 & 0.06 & 0.06 & 0.06 \\ 
3 & 26 & 0.06 & 0.06 & 0.07 & 0.07 & 0.06 & 0.06 & 0.06 & 0.06 & 0.06 \\ 
3 & 31 & 0.03 & 0.03 & 0.04 & 0.04 & 0.03 & 0.03 & 0.03 & 0.03 & 0.03 \\ 
3 & 36 & 0.04 & 0.04 & 0.05 & 0.05 & 0.04 & 0.04 & 0.04 & 0.04 & 0.04 \\ 
3 & 41 & 0.03 & 0.03 & 0.04 & 0.04 & 0.03 & 0.03 & 0.03 & 0.03 & 0.03 \\ 
3 & 46 & 0.04 & 0.04 & 0.04 & 0.04 & 0.04 & 0.04 & 0.04 & 0.04 & 0.04 \\ 
3 & 51 & 0.04 & 0.04 & 0.04 & 0.04 & 0.04 & 0.04 & 0.04 & 0.04 & 0.04 \\ 
6 & 6 & 0.05 & 0.05 & 0.05 & 0.05 & 0.05 & 0.05 & 0.05 & 0.05 & 0.05 \\ 
6 & 11 & 0.04 & 0.04 & 0.04 & 0.04 & 0.04 & 0.04 & 0.04 & 0.04 & 0.04 \\ 
6 & 16 & 0.03 & 0.03 & 0.03 & 0.03 & 0.03 & 0.03 & 0.03 & 0.03 & 0.03 \\ 
6 & 21 & 0.05 & 0.05 & 0.05 & 0.05 & 0.05 & 0.05 & 0.05 & 0.05 & 0.05 \\ 
6 & 26 & 0.05 & 0.05 & 0.06 & 0.06 & 0.05 & 0.05 & 0.05 & 0.05 & 0.05 \\ 
6 & 31 & 0.03 & 0.03 & 0.03 & 0.03 & 0.03 & 0.03 & 0.03 & 0.03 & 0.03 \\
\hline 
\end{tabular}
\label{table:bektasCMILP1}
\end{footnotesize}
\end{table*}

\begin{table*}[htbp]
\center
\begin{footnotesize}
\caption{The average computational times (in seconds) for the combinations of CMILP+TSP algorithms proposed by \cite{bektacs2014} on 80 instances}
\begin{tabular}{ccccccccccc}
\hline
 &  & CMILP+ & CMILP+ & CMILP+ & CMILP+ & CMILP+ & CMILP+ & CMILP+ & CMILP+ & CMILP+DL+ \\ 
n & r & DL+3CLQ & DL+NR &  DL+L3 & DL+2P & DL+R & R+2P & NR+2P & NR+R+2P &  NR+R+2P \\ 
\hline
3 & 6 & 0.57 & 0.58 & 0.44 & 0.57 & 0.57 & 0.51 & 0.52 & 0.70 & 0.74 \\ 
3 & 11 & 2.21 & 2.73 & 1.43 & 2.66 & 2.68 & 4.59 & 4.32 & 6.31 & 6.32 \\ 
3 & 16 & 6.66 & 8.67 & 3.63 & 8.78 & 9.09 & 16.41 & 15.73 & 23.18 & 23.98 \\ 
3 & 21 & 15.53 & 20.39 & 8.39 & 22.86 & 20.04 & 42.88 & 41.30 & 65.15 & 65.59 \\ 
3 & 26 & 31.04 & 40.70 & 16.93 & 50.69 & 41.09 & 98.33 & 94.67 & 154.04 & 155.82 \\ 
3 & 31 & 56.39 & 73.05 & 32.15 & 98.04 & 73.78 & 198.39 & 198.67 & 308.93 & 311.91 \\ 
3 & 36 & 97.01 & 129.09 & 58.51 & 204.77 & 125.29 & 397.43 & 385.07 & 623.18 & 632.16 \\ 
3 & 41 & 165.45 & 205.95 & 104.06 & 374.08 & 208.42 & 724.29 & 711.80 & 1170.70 & 1189.32 \\ 
3 & 46 & 269.28 & 326.19 & 178.89 & 594.52 & 330.80 & 1331.37 & 2306.22 & 2178.78 & 2568.61 \\ 
3 & 51 & 429.35 & 506.52 & 298.81 & 1087.05 & 502.58 & 2290.12 & 2282.70 & 7378.42 & 8390.89 \\ 
6 & 6 & 0.56 & 0.73 & 0.30 & 0.51 & 0.81 & 0.94 & 0.90 & 1.30 & 1.33 \\ 
6 & 11 & 3.30 & 5.27 & 1.75 & 4.28 & 4.91 & 8.78 & 8.35 & 12.33 & 13.08 \\ 
6 & 16 & 11.56 & 20.08 & 7.05 & 16.44 & 16.64 & 34.11 & 32.62 & 49.62 & 50.77 \\ 
6 & 21 & 29.19 & 40.24 & 17.24 & 45.71 & 41.31 & 94.47 & 91.77 & 146.95 & 150.05 \\ 
6 & 26 & 63.14 & 83.91 & 36.51 & 117.37 & 86.66 & 243.64 & 232.01 & 380.74 & 382.88 \\ 
6 & 31 & 126.60 & 161.97 & 77.62 & 256.15 & 165.80 & 549.37 & 535.47 & 866.40 & 892.54 \\ 
\hline
\end{tabular}
\label{table:bektasCMILP2}
\end{footnotesize}
\end{table*}

\begin{table*}[htbp]
\center
\begin{footnotesize}
 \caption{The average optimality gap(\%) for the combinations of SP+TSP algorithms proposed by \cite{bektacs2014} on 80 instances}

\begin{tabular}{ccccccccccc}
\hline
\multicolumn{1}{l}{} & \multicolumn{1}{l}{} & SP+ & SP+ & SP+ & SP+ & SP+ & SP+ & SP+ & SP+ & SP+DL+ \\ 
n & r & DL+3CLQ & DL+NR &  DL+L3 & DL+2P & DL+R & R+2P & NR+2P & NR+R+2P &  NR+R+2P \\ 
\hline
3 & 6 & 0.06 & 0.06 & 0.06 & 0.06 & 0.06 & 0.06 & 0.06 & 0.06 & 0.06 \\ 
3 & 11 & 0.05 & 0.05 & 0.05 & 0.05 & 0.05 & 0.05 & 0.05 & 0.05 & 0.05 \\ 
3 & 16 & 0.04 & 0.04 & 0.04 & 0.04 & 0.04 & 0.04 & 0.04 & 0.04 & 0.04 \\ 
3 & 21 & 0.06 & 0.06 & 0.06 & 0.06 & 0.06 & 0.06 & 0.06 & 0.06 & 0.06 \\ 
3 & 26 & 0.07 & 0.07 & 0.08 & 0.07 & 0.07 & 0.07 & 0.07 & 0.07 & 0.07 \\ 
3 & 31 & 0.03 & 0.03 & 0.04 & 0.04 & 0.03 & 0.03 & 0.03 & 0.03 & 0.03 \\ 
3 & 36 & 0.04 & 0.04 & 0.05 & 0.05 & 0.04 & 0.04 & 0.04 & 0.04 & 0.04 \\ 
3 & 41 & 0.03 & 0.03 & 0.04 & 0.04 & 0.03 & 0.03 & 0.03 & 0.03 & 0.03 \\ 
3 & 46 & 0.04 & 0.04 & 0.04 & 0.04 & 0.04 & 0.04 & 0.04 & 0.04 & 0.04 \\ 
3 & 51 & 0.04 & 0.04 & 0.04 & 0.04 & 0.04 & 0.04 & 0.04 & 0.04 & 0.04 \\ 
6 & 6 & 0.05 & 0.05 & 0.05 & 0.05 & 0.05 & 0.05 & 0.05 & 0.05 & 0.05 \\ 
6 & 11 & 0.04 & 0.04 & 0.04 & 0.04 & 0.04 & 0.04 & 0.04 & 0.04 & 0.04 \\ 
6 & 16 & 0.04 & 0.04 & 0.04 & 0.04 & 0.04 & 0.04 & 0.04 & 0.04 & 0.04 \\ 
6 & 21 & 0.05 & 0.05 & 0.05 & 0.05 & 0.05 & 0.05 & 0.05 & 0.05 & 0.05 \\ 
6 & 26 & 0.06 & 0.06 & 0.06 & 0.06 & 0.06 & 0.06 & 0.06 & 0.06 & 0.06 \\ 
6 & 31 & 0.03 & 0.03 & 0.03 & 0.03 & 0.03 & 0.03 & 0.03 & 0.03 & 0.03 \\ 
\hline
\end{tabular}
\label{table:bektasSP2}
\end{footnotesize}
\end{table*}


\begin{table*}[h]
\center
\begin{footnotesize}
\caption{The average computational times (in seconds) for the combinations of SP+TSP algorithms proposed by \cite{bektacs2014} on 80 instances}
\begin{tabular}{ccccccccccc}
\hline
\multicolumn{1}{l}{} & \multicolumn{1}{l}{} & SP+ & SP+ & SP+ & SP+ & SP+ & SP+ & SP+ & SP+ & SP+DL+ \\ 
n & r & DL+3CLQ & DL+NR &  DL+L3 & DL+2P & DL+R & R+2P & NR+2P & NR+R+2P &  NR+R+2P \\ 
\hline
3 & 6 & 0.36 & 0.34 & 0.21 & 0.34 & 0.34 & 0.52 & 0.57 & 0.73 & 0.78 \\ 
3 & 11 & 1.83 & 2.36 & 1.04 & 2.32 & 2.40 & 4.41 & 4.21 & 6.11 & 6.52 \\ 
3 & 16 & 6.31 & 8.48 & 3.47 & 8.76 & 8.42 & 16.69 & 15.84 & 23.79 & 23.87 \\ 
3 & 21 & 16.26 & 21.21 & 9.48 & 23.50 & 21.26 & 44.73 & 42.05 & 63.84 & 64.88 \\ 
3 & 26 & 34.15 & 44.64 & 20.70 & 53.14 & 45.06 & 99.55 & 95.77 & 149.88 & 148.41 \\ 
3 & 31 & 65.81 & 82.44 & 43.16 & 107.32 & 83.84 & 209.11 & 199.04 & 315.52 & 323.29 \\ 
3 & 36 & 117.06 & 142.98 & 78.11 & 207.95 & 146.99 & 411.47 & 399.74 & 636.82 & 642.91 \\ 
3 & 41 & 203.06 & 240.48 & 140.42 & 374.23 & 246.66 & 778.39 & 917.59 & 1204.61 & 1221.43 \\ 
3 & 46 & 328.83 & 388.57 & 242.08 & 674.43 & 396.26 & 1390.94 & 1419.97 & 2291.09 & 2260.42 \\ 
3 & 51 & 528.08 & 614.95 & 396.43 & 1229.09 & 626.53 & 2442.66 & 2578.34 & 3980.35 & 4007.79 \\ 
6 & 6 & 0.59 & 0.65 & 0.38 & 0.61 & 0.68 & 1.01 & 1.01 & 1.42 & 1.43 \\ 
6 & 11 & 4.02 & 5.05 & 2.28 & 4.91 & 5.16 & 9.18 & 9.04 & 13.16 & 13.34 \\ 
6 & 16 & 15.74 & 18.17 & 8.81 & 19.13 & 18.98 & 35.88 & 35.64 & 52.29 & 52.90 \\ 
6 & 21 & 38.83 & 63.85 & 24.93 & 56.57 & 50.07 & 102.75 & 102.27 & 155.33 & 155.32 \\ 
6 & 26 & 88.07 & 106.35 & 59.02 & 134.84 & 111.04 & 264.43 & 258.25 & 401.40 & 400.82 \\ 
6 & 31 & 179.65 & 211.39 & 127.40 & 308.35 & 219.98 & 597.99 & 586.00 & 932.82 & 934.03 \\ 
\hline
\end{tabular}
\label{table:bektasSP1}
\end{footnotesize}
\end{table*}

\begin{table*}[htbp]
\center
\begin{footnotesize}
\caption{Notation of formulations}
\begin{tabular}{l|llllr}
\hline
CMILP & \multicolumn{ 4}{l}{Classical mixed integer linear programming} &  \\ 
\hline
SP  & \multicolumn{ 3}{l}{Shortest path flow formulation } & \multicolumn{1}{r}{} &  \\ 
\hline
MTZ & \multicolumn{ 4}{l}{TSP formulation proposed by \cite{miller1960}} &  \\ 
\hline
DL & \multicolumn{ 4}{l}{TSP formulation proposed by \cite{desrochers1991}} &  \\ 
\hline
SC & \multicolumn{ 5}{l}{Single commodity flow formulation proposed by \cite{gavish1978} }  \\ 
\hline
2C & \multicolumn{ 5}{l}{Two commodity flow formulation proposed by \cite{finke1984} }   \\  
\hline
MTZ+2CLQ & \multicolumn{ 5}{l}{} \\ 
3CLQ & \multicolumn{ 5}{l}{} \\ 
NR & \multicolumn{ 5}{r}{The generalization of MTZ formulation proposed by \cite{bektacs2014}} \\  
L3 & \multicolumn{ 5}{r}{} \\ 
2P & \multicolumn{ 5}{r}{} \\ 
R & \multicolumn{ 5}{r}{} \\ 
\hline
\end{tabular}
\label{}
\end{footnotesize}
\end{table*}

\end{document}